\documentclass[12pt]{article}
\usepackage{amsfonts}
\usepackage{amssymb,amsmath,amsthm,latexsym}
\usepackage[numbers,compress]{natbib}
\usepackage{amsmath}
\usepackage{cases}
\usepackage{indentfirst}
\usepackage{mathrsfs}
\usepackage{amsfonts}
\usepackage{amsmath}
\usepackage{amsthm}
\usepackage[usenames]{color}
\usepackage{footnote}
\usepackage{longtable}
\usepackage{hyperref}
\usepackage[all]{xy}
\usepackage{amsmath,amscd}
\usepackage{multirow} 
\usepackage{booktabs}
\usepackage{clrscode}
\usepackage{listings}

\allowdisplaybreaks[4]
\newtheorem{theorem}{Theorem}[section]

\newtheorem{corollary}[theorem]{Corollary}

\newtheorem{define}[theorem]{Definition}

\newcommand{\pf} { {\rm \noindent{\bf Proof.}} }

\numberwithin{equation}{section}

\begin{document}

\title{ Permutation polynomials and complete permutation polynomials over $\mathbb{F}_{q^{3}}$ }
\author{Yanping Wang$^{1,4}$, Weiguo Zhang$^{1,2,}${\thanks{Corresponding author.  Email addresses:~zwg@xidian.edu.cn} }, Daniele Bartoli$^3$, Qiang Wang$^4$ \\
\vspace*{0.0cm}\\
{\small $^1$ISN Laboratory, Xidian University, Xi'an, 710071, China}\\
{\small $^2$State Key Laboratory of Cryptology, P.O. Box 5159, Beijing 100878, China}\\
{\small $^3$Department of Mathematics and Computer Science, University of Perugia, Perugia, 06123, Italy }\\
{\small $^4$School of Mathematics and Statistics, Carleton University, Ottawa, ON K1S 5B6, Canada}\\
}
\date{}
\maketitle
\begin{abstract}
Motivated by many recent constructions of permutation polynomials over $\mathbb{F}_{q^2}$, we study
permutation polynomials over $\mathbb{F}_{q^3}$ in terms of their coefficients.
Based on the multivariate method and resultant elimination, we construct several new classes of sparse permutation polynomials over $\mathbb{F}_{q^3}$, $q=p^{k}$, $p\geq3$.
Some of them are complete mappings.
\end{abstract}




\section{Introduction}
Let $q$ be a power of prime $p$ and $\mathbb{F}_q$ denote the finite field with $q$ elements. Define $\mathbb{F}^{*}_q$ to be the multiplicative group of $\mathbb{F}_q$.
A polynomial $f(x)\in \mathbb{F}_q[x]$ is called a \emph{permutation polynomial} over $\mathbb{F}_q$ if the associated polynomial function $f: c \mapsto f(c)$ from $\mathbb{F}_q$ into $\mathbb{F}_q$ is a permutation of $\mathbb{F}_q$ \cite{LN97}. A permutation polynomial $f(x)\in \mathbb{F}_q[x]$ is a \emph{complete permutation polynomial} over $\mathbb{F}_q$ if $f(x)+x$ permutes $\mathbb{F}_q$ as well. In general, for an $\epsilon\in \mathbb{F}_q^*$ a permutation polynomial $f(x)$ over $\mathbb{F}_q$ is called a $\epsilon$-complete permutation polynomial if $f(x)+\epsilon x$ is also a permutation polynomial. Note that if $f(x)$ is an $\epsilon$-complete permutation polynomial then $\epsilon^{-1}f(x)$ is a complete permutation polynomial.

Finding new permutation polynomials and complete permutation polynomials is of great interest in both theoretical and applied aspects. Many constructions of permutation polynomials appeared in the recent years; see for instance  \cite{AGW2011,BLZV2015,Hou2015,TW2017,WL2015,XLZH2018}.  The reader may refer to {\rm\cite[Chapter $7$]{LN97}}, {\rm\cite[Chapter $8$]{GD2013}},   \cite{houx2015} and references therein for more information.

Permutation polynomials with few terms are interesting for their simple algebraic forms and have wide applications in coding theory {\rm\cite{Laigle2007}}, combinatorial designs {\rm \cite{LN97}}, and cryptography {\rm \cite{G1993}}. Recently,  the multivariate method (a key tool in the proof of Niho's conjecture {\rm\cite{hans1999}}) has been used to construct permutation polynomials with few terms.  Dobbertin developed the multivariate method {\rm\cite{DOBBERTIN2002}} to confirm the permutation property of certain types of polynomials over $\mathbb{F}_{2^{n}}$.
Later on, Ding et al.  {\rm\cite{CLQ2015}} explored the multivariate method to construct several classes of permutation trinomials over finite fields with even characteristic.
Motivated by what Ding et al. did for permutation trinomials in {\rm\cite{CLQ2015}}, Li et al. {\rm\cite{LQC2017}} proposed several classes of permutation trinomials over $\mathbb{F}_{2^n}$. Three more classes of permutation trinomials were given by Ma et al. {\rm\cite{JTTG2017}} as well. Recently Bartoli and Zini {\rm\cite{BZ2018}} determined all permutation trinomials of the form $x^{2p^s+r}+x^{p^s+r}+\lambda x^r$ over $\mathbb{F}_{p^{t}}$ when $(2p^s+r)^4<p^t$. Wang et al. {\rm\cite{WZZ2017,WZZ2018}} presented several classes of permutation trinomials over $\mathbb{F}_{q^{n}}$ with characteristic $2$ and $3$. Furthermore, Bartoli {\rm\cite{Bartoli2018}} characterized four classes of permutation trinomials over $\mathbb{F}_{q^3}$ in terms of their coefficients in $\mathbb{F}_{q}$, $q=p^k$ and $p>3$.

Other types of approaches can be seen as applications of the AGW criterion {\rm\cite{AGW2011}}. In particular, polynomials of type
\begin{eqnarray}\label{eq:form}
x^{r}f(x^{(q-1)/d})
\end{eqnarray}
 over the finite field $\mathbb{F}_{q}$ has been studied  earlier  by Wan and Lidl {\rm\cite{WR1991}}, Park and Lee {\rm\cite{PL2001}},  Akbary and Wang {\rm\cite{AW2007}},  Wang {\rm\cite{WQ2007}}, and Zieve  {\rm\cite{MEZ2009}}. A polynomial of the form (\ref{eq:form}) permutes $\mathbb{F}_{q}$ if and only if $\gcd(r, (q-1)/d)=1$ and $x^{r}f(x)^{(q-1)/d}$ permutes the set $\mu_{d}$ of the $d$-th roots of unity in $\mathbb{F}_{q}$. Many  classes of permutation polynomials are characterized using this criterion.
For example, using this approach, permutation trinomials over $\mathbb{F}_{q^{2}}$ from Niho exponents
$$ ax^{s(2^{m}-1)+1} + bx^{t(2^{m}-1)+1} + x,$$
have been characterized; see {\rm\cite{hou2014,Hou2015,RR2016,ZHF2017,LT2017,DZ2018,DYDM2016,BG2018,TZLH2018}}.
In particular,  Hou {\rm\cite{hou2014,Hou2015}} completely determined the permutation behaviors of trinomials with $(s, t)=(1, 2)$ over $\mathbb{F}_{q^{2}}$ by discussing all possible cases of the coefficients $a$ and $b$.
Tu et al. {\rm\cite{TZLH2018}} characterized a class of new permutation trinomials with with $(s, t)=(q, 2)$, $q=2^m$, in terms of the coefficients $a$ and $b$ over $\mathbb{F}_{q^{2}}$.  They proved the sufficiency of the conditions  and conjectured their necessity.  Then Bartoli \cite{Bartoli:18} proved the necessity using low degree algebraic curves and computational packages such as MAGMA. Hou \cite{Hou:18} found a way to prove both directions at the same time.

In this paper, we study permutation polynomials over $\mathbb{F}_{q^3}$ of types as shown in Table \ref{Table1}.
\begin{table}
\caption{Types of permutation polynomials}\label{Table1}
\begin{center}
\begin{tabular}{|l|l|}
\hline
Polynomials & References\\
\hline\hline
$x^{q^2+q-1} + Ax$&Theorem \ref{th:pb}\\
\hline
$x^{q^2-q+1} + Ax^{q^3-q^2+q} + Bx$&Theorem \ref{th:pt1}\\
\hline
$x^{q^2+q-1} + Ax^{q^2} + Cx$&Theorem \ref{th:pt2}\\
\hline
$x^{q^2+q-1} + Bx^{q} + Cx$&Theorem \ref{th:pt3}\\
\hline
$x^{q^2+q-1} + Ax^{q^2} + Bx^{q} + Cx$&Theorem \ref{th:pq1}\\
\hline
$x^{q^2+q-1} + Ax^{q^2-q+1} + Bx^{q^2} + Cx$&Theorem \ref{th:pq3}\\
\hline
$x^{q^2+q-1} + Ax^{q^3-q^2+q} + Bx^{q} + Cx$&Theorem \ref{th:pq4}\\
\hline
$x^{q^2+q-1} + Ax^{q^2-q+1} + Bx^{q} + Cx$&Theorem \ref{th:pq2}\\
\hline
$x^{q^2+q-1} + Ax^{q^2-q+1} + Bx^{q^2} + Cx^{q} + Dx$&Theorem \ref{th:pp1}\\
\hline
\end{tabular}
\end{center}
\end{table}
The paper is organized as follows. Section \ref{pre} gives some preliminaries on the resultant between two polynomials. In Section \ref{three}, by using the multivariate method, a class of complete permutation binomials, three classes of permutation trinomials, four classes of permutation quadrinomials and one class of permutation pentanomials over $\mathbb{F}_{q^3}$ are costructed in terms of their coefficients.
The conclusion is presented in Section \ref{conclu}. Some necessary Magma programs are included in the Appendix \ref{App}.

\section{Preliminaries }\label{pre}
In some of our proofs we will need to investigate the solutions of a system of polynomial equations. In these situation an important tool we use is the resultant of two polynomials.
We recall here some basic facts about the resultant.

\begin{define}({\rm \cite{LN97}}, p.36) Let $f(x)=a_{0}x^{n}+a_{1}x^{n-1}+\cdots +a_{n}\in \mathbb{F}_q[x]$ and $g(x)=b_{0}x^{m}+b_{1}x^{m-1}+\cdots+ b_{m}\in \mathbb{F}_q[x]$ be
two polynomials of formal degree $n$ respectively $m$ with $n, m \in \mathbb{N}$. Then the resultant $R(f, g)$ of the two polynomials is defined by the determinant
\begin{eqnarray*}
R(f, g)=
\left|\begin{array}{cccccccc}
   a_{0} &    a_{1}    & \cdots &   a_{n}  & 0  &    &\cdots & 0   \\
    0  &  a_{0} &    a_{1}    & \cdots &   a_{n}  & 0  & \cdots & 0 \\
    \vdots  &    &        &      &     &      &          &  \vdots    \\
    0  &  \cdots &   0   &  a_{0} &    a_{1}  &     &  \cdots &   a_{n} \\
    b_{0} &    b_{1}    & \cdots &      & b_{m}  & 0  & \cdots &  0     \\
      0   &    b_{0} &   b_{1}  & \cdots &    &  b_{m}  & \cdots &  0 \\
    \vdots &          &      &     &      &     &    &\vdots       \\
    0  &   \cdots   & 0 &  b_{0} &    b_{1} &    & \cdots &  b_{m}  \\
\end{array}\right|
& \begin{array}{l}
\left.\rule{0mm}{10.5mm}\right\}$m$~rows \\
\\\left.\rule{0mm}{10.5mm}\right\}$n$~ rows
\end{array}\\[0pt]
\end{eqnarray*}
of order $m+n$.
\end{define}
If $deg(f) = n$ (i.e., if $a_{0}\neq 0$) and $f(x)=a_{0}(x -\alpha_{1})(x -\alpha_{2})\cdots (x -\alpha_{n})$ in
the splitting field of $f$ over $\mathbb{F}_q$, then $R(f, g)$ is also given by the formula
\begin{eqnarray*}
R(f, g)=a_{0}^{m}\prod_{i=1}^n g(\alpha_{i}).
\end{eqnarray*}
In this case, we obviously have $R(f, g) = 0$ if and only if $f$ and $g$ have a common root, which is the same as saying that $f$ and $g$ have a common
divisor in $\mathbb{F}_q[x]$ of positive degree.

For  two polynomials $F(x,y),\, G(x,y)\in \mathbb{F}_q[x,y]$ of positive degree in $y$   the resultant $Res(F,G,y)$ of $F$ and $G$ with respect to $y$  is the resultant of $F$ and $G$ when considered as polynomials in the single variable $y$ (that is, as elements in $R[y]$ with $R=\mathbb{F}_q[x]$). In this case $Res(F,G,y)\in \mathbb{F}_q[x]$ is in the ideal generated by $F$ and $G$, and therefore any pair $(a,b)$ with $F(a,b)=G(a,b)=0$ is such that $Res(F,G,y)(a)=0$; see e.g. \cite[Prop 3.6.1]{COX}.

Let $d$ be a divisor of $q^3-1$. We denote by  $\mu_{d}=\{x \in \mathbb{F}_{q^3} \ : \ x^d =1\}$ the set of $d$-th roots of unity in $\mathbb{F}_{q^3}$.

\section{The classes of permutation polynomials and complete permutation polynomials over $\mathbb{F}_{q^3}$}\label{three}
In the section, we characterize a class of complete permutation binomials over $\mathbb{F}_{q^3}$, three classes of permutation trinomials, and four classes of permutation quadrinomials and one class of permutation pentanomials over $\mathbb{F}_{q^3}$ in terms of their coefficients over $\mathbb{F}_{q}$. The multivariate method as a very useful tool is applied to prove these new results.

\subsection{A class of complete permutation binomials }
One class of complete permutation binomials is presented over $\mathbb{F}_{q^3}$ in the subsection.

\begin{theorem}\label{th:pb}
Let $A\in\mathbb{F}_{q^3}$ with $A^{q^2+q+1}=-1$, $(A+1)^{q^2+q+1}=-1$, and $A, A+1\notin \mu_{\frac{q^2+q+1}{3}}$. Then
\begin{eqnarray}\label{eq:bino}
f(x)= x^{q^2+q-1} + Ax
\end{eqnarray}
is a complete permutation binomial over $\mathbb{F}_{q^3}$.
\end{theorem}

\pf We will prove that for each $a \in\mathbb{F}_{q^3}$, the equation
\begin{eqnarray}\label{eq:b1.1}
x^{q^2+q-1} + Ax = a
\end{eqnarray}
has at most a solution in $\mathbb{F}_{q^3}$.  We demonstrate $x=0$ if and only if $a=0$. Supposing $x\neq0$ is one solution of equation
\begin{eqnarray}\label{eq:b1.2}
x^{q^2+q-1} + Ax =0,
\end{eqnarray}
which is necessary to prove that
\begin{eqnarray}\label{eq:b1.3}
x^{q^2+q-2} + A =0
\end{eqnarray}
has no solution.
Setting $u=x^{q-1}$ implies $u^{q^2+q+1}=1$. From Eq. (\ref{eq:b1.3}) we have $u^{q} = -\frac{A}{u^2}$ and $u^{q^2} = -\frac{A^{q-2}}{u^4}$.
Substituting them into $u^{q^2+q+1}=1$, we get
\begin{eqnarray*}
u^{3} = -\frac{1}{A^{q-1}}.
\end{eqnarray*}
Since $\gcd(q+2, 3)=3$, substituting $u^{3} = -\frac{1}{A^{q-1}}$ into the equation $u^{q+2}+A=0$, we obtain
$$ A^{\frac{q^2+q+1}{3}}=1, $$
which contradicts $A$ is not the $\frac{q^2+q+1}{3}$ roots of unity in $\mathbb{F}_{q^3}$.
So $f(x)=0$ if and only if $x=0$.

If $a\neq0$, we show that Eq. (\ref{eq:b1.1}) has one nonzero solution. Let $y=x^{3^{k}}$, $z=y^{3^{k}}$, $b=a^{3^{k}}$ and $c=b^{3^{k}}$, then we obtain the system of equations
\begin{numcases}{}
yz + Ax^{2} -ax =0,  \label{eq:b1.11}\\
zx + A^{q}y^{2} -by =0,  \label{eq:b1.11'}\\
xy + A^{q^2}z^{2} -cz =0.  \label{eq:b1.11''}
\end{numcases}
Let $B=A^{q}$, $C=A^{q^2}$, then $C^{q}=A$ and through a series of computations of the resultant with $ABC=-1$ (see \ref{App}), we have
\begin{eqnarray*}
&& B^{3}x^{7}(Ax-a)\big[(- a^3B^3C^3 - 2abcAB^2C^2 + abcBC + b^3AC^2 + c^3AB^2)x \nonumber \\
&+&  a^2bcB^2C^2 - ab^3C^2 - ac^3B^2 + b^2c^2 \big]=0.
\end{eqnarray*}
Since $a\neq0$, $x=0$ is not a solution. If $x= \frac{a}{A}$, then substituting $x=\frac{a}{A}$ into Eq. (\ref{eq:b1.1}) gets $x=0$, which is impossible,
therefore the equation has at most one solution.

Similar computations show that  the equation $x^{q^2+q-1} + (A+1)x = a$, where $a \in\mathbb{F}_{q^3}$, has at most a solution in $\mathbb{F}_{q^3}$ if $(A+1)^{q^2+q+1}=-1$ and $A+1\notin \mu_{\frac{q^2+q+1}{3}}$ and the claim follows.
\qed

\begin{corollary}\label{cor:pcb}
Let $q=p^k$, $k\in \mathbb{N}$. The number of $A$ is $\frac{2}{3}(q^2+q+1)$ such that
\begin{eqnarray}\label{eq:bino}
f(x)= x^{q^2+q-1} + Ax
\end{eqnarray}
is a complete permutation binomial over $\mathbb{F}_{q^3}$, where
$A=\theta^{\frac{m(q-1)}{2}}$ with $\theta\in\mathbb{F}^{*}_{q^3}$  is a primitive element of $\mathbb{F}_{q^3}$ and $m\equiv  \pm 1 ({\rm mod}\ 6)$.
\end{corollary}

\subsection{Three classes of permutation trinomials }
Three classes of permutation trinomials of the form $f(x)=x^{d_{1}} + Ax^{d_{2}} + Bx \in\mathbb{F}_{q}[x]$ are given in the subsection.

\begin{theorem}\label{th:pt1}
Let $A, B\in\mathbb{F}_{q}$ with $A, B\neq0$, $AB-1\neq0$, $AB+2\neq0$, $A + B + 1\neq0$, $A^2 - AB - A + B^2 - B + 1\neq0$ and $A^{3}+AB+1=0$. Then
\begin{eqnarray}\label{eq:pt1.1}
f(x)= x^{q^2-q+1} + Ax^{q^3-q^2+q} + Bx
\end{eqnarray}
is a permutation trinomial over $\mathbb{F}_{q^3}$.
\end{theorem}

\pf We show that for each $a \in\mathbb{F}_{q^3}$, the equation
\begin{eqnarray}\label{eq:t2.1}
x^{q^2-q+1} + Ax^{q^3-q^2+q} + Bx = a
\end{eqnarray}
has at most a solution in $\mathbb{F}_{q^3}$.  We prove $x=0$ if and only if $a=0$. Suppose that $x\neq0$ is one solution of equation
\begin{eqnarray}\label{eq:t2.2}
x^{q^2-q+1} + Ax^{q^3-q^2+q} + Bx =0.
\end{eqnarray}
Setting $u=x^{q-1}$ implies $u^{q^2+q+1}=1$. It is necessary to prove that
\begin{eqnarray}\label{eq:t2.3}
u^{q} + Au^{q^2+1} + B =0
\end{eqnarray}
has no solution.
Raising Eq. (\ref{eq:t2.3}) to the $q$-th power results in
\begin{eqnarray}\label{eq:t2.4}
u^{q^2} + Au^{q+1} + B =0
\end{eqnarray}
Combining Eqs. (\ref{eq:t2.3}) and (\ref{eq:t2.4}) obtains
\begin{eqnarray}\label{eq:t2.5}
A^{2}u^{q+2} - u^q + ABu - B =0,
\end{eqnarray}
therefore we deduce $u^{q} = \frac{B(1-Au)}{A^{2}u^2-1}$ and
$u^{q^2} =\frac{B(A^{2}u^2-1)^{2}-AB^2(1-Au)(A^{2}u^2-1)}{(A^{2}Bu-AB)^{2}-(A^{2}u^2-1)^{2}}$.
If $A^{2}u^2-1=0$, then $u=\frac{1}{A}$ or $u=-\frac{1}{A}$.
Again substituting $u=\frac{1}{A}$ into Eq. (\ref{eq:t2.3}) results in $AB+2=0$, which is a contradiction. Similarly, we plug $u=-\frac{1}{A}$ into Eq. (\ref{eq:t2.3}) obtains $B=0$, which is impossible.
If $(A^{2}Bu-AB)^{2}-(A^{2}u^2-1)^{2}=0$, then we obtain $A^{2}u^2+A^{2}Bu-AB-1=0$ or $A^{2}u^2-A^{2}Bu+AB-1=0$.
We only discuss the first equation, the other is similarly discussed. Combining Eq. (\ref{eq:t2.5}) and the first equation derives
$$B(Au-1)(Au^q-1)=0.$$
Since $B\neq 0$ and $u\neq\frac{1}{A}$, thereby we have $u^q=\frac{1}{A}$ which contradicts $u\neq\frac{1}{A}$.

By $u^{q^2+q+1}=1$ and $A^{3}+AB+1=0$ we have
\begin{eqnarray}\label{eq:t2.6}
(Au - 1)^3(Au + 1)(Au + AB + 1)((A - B^2)u - AB + 1)= 0,
\end{eqnarray}
since $AB+2\neq0$ and $B\neq0$ means that $u\neq\frac{1}{A}$ and $u\neq-\frac{1}{A}$. Therefore we have
\begin{eqnarray}\label{eq:t2.7}
Au + AB + 1=0, ~~~~~~~~(A - B^2)u - AB + 1= 0.
\end{eqnarray}
If $Au+AB+1=0$, then substituting $u=-\frac{AB+1}{A}$ into Eq. (\ref{eq:t2.3}) derives $AB(AB+2)=0$, which contradicts
$A, B\neq0$ and $AB+2\neq0$. Since $AB\neq1$, hence we obtain $u=\frac{AB - 1}{A - B^2}$. Substituting $u=\frac{AB - 1}{A - B^2}$ into Eq. (\ref{eq:t2.3}) leads to $B^2(A + B + 1)(A^2 - AB - A + B^2 - B + 1)=0$. Since $A + B + 1\neq0$ and $A^2 - AB - A + B^2 - B + 1\neq0$, it is a contradiction, hence $f(x)=0$ has one unique solution $x=0$.

If $a\neq0$, then we show that Eq. (\ref{eq:t2.1}) has one nonzero solution. Let $y=x^{q}$, $z=y^{q}$, $b=a^{q}$ and $c=b^{q}$, then we obtain the system of equations
\begin{numcases}{}
xz^2 + Axy^2 + Bxyz -ayz =0,  \label{eq:t2.11}\\
yx^2 + Ayz^2 + Bxyz -bzx =0,  \label{eq:t2.11'}\\
zy^2 + Azx^2 + Bxyz -cxy =0.  \label{eq:t2.11''}
\end{numcases}
Eliminating $z$ by (\ref{eq:t2.11}) and (\ref{eq:t2.11''}), we have
\begin{eqnarray*}
f_{1} \triangleq \xi_{1}y^4 + \xi_{2}y^3 + \xi_{3}y^2 + \xi_{4}y + \xi_{5}=0,
\end{eqnarray*}
where $\xi_{i}:= \xi_{i}(x, A, B, a, b, c)$, $i= 1, 2, 3, 4, 5$.

Furthermore, by Eqs. (\ref{eq:t2.11'}) and (\ref{eq:t2.11''}) we get
\begin{eqnarray*}
f_{2} \triangleq \eta_{1}y^4 + \eta_{2}y^3 + \eta_{3}y^2 + \eta_{4}y + \eta_{5}=0,
\end{eqnarray*}
where $\eta_{i}:= \eta_{i}(x, A, B, a, b, c)$, $i= 1, 2, 3, 4, 5$.

Computing the resultant of $f_{1}$ and $f_{2}$ with respect for $y$ and recalling that $A^{3}= -AB-1$, we derive
\begin{eqnarray*}
c^8x^{4}(\alpha x +\beta)=0,
\end{eqnarray*}
where $\alpha := \alpha(A, B, a, b, c)$, $\beta := \beta(A, B, a, b, c)$. Since $a\neq0$ implies $c\neq0$ , thus $x=0$ is not a solution of the above equation and the equation has at most a solution $x=-\frac{\beta}{\alpha}$.
Therefore Eq. (\ref{eq:t2.1}) has at most a solution and we complete the proof.
\qed

In the following we present two classes of permutation trinomials with the similar form which have been discussed in {\rm\cite{Bartoli2018}}, however, the two classes of polynomials that we will discuss next require different restrictions on coefficients.

\begin{theorem}\label{th:pt2}
Let $A, C\in\mathbb{F}_{q}$ with $C^2-C+1=0$, $A^{3}\neq-1$ and $m_{A,C}(t)=Ct^3 + A(C+1)t^2 + A^{2}t - C\in \mathbb{F}_{q}[t]$ has no roots in $\mu_{q^2+q+1}$. Then
\begin{displaymath}
f(x)= x^{q^2+q-1} + Ax^{q^2} + Cx
\end{displaymath}
is a permutation trinomial over $\mathbb{F}_{q^3}$.
\end{theorem}

\pf As in the previous theorems, we will show that for each $a \in\mathbb{F}_{q^3}$, the equation
\begin{eqnarray}\label{eq:pt2.1}
x^{q^2+q-1} + Ax^{q^2} + Cx = a
\end{eqnarray}
has at most a solution in $\mathbb{F}_{q^3}$.  We prove $x=0$ if and only if $a=0$. Supposing $x\neq0$ is one solution of equation
\begin{eqnarray}\label{eq:pt2.2}
x^{q^2+q-1} + Ax^{q^2} + Cx =0.
\end{eqnarray}
Setting $u=x^{q-1}$ implies $u^{q^2+q+1}=1$. It is necessary to prove that
\begin{eqnarray}\label{eq:pt2.3}
u^{q+2} + Au^{q+1} + C =0
\end{eqnarray}
has no solution.
From Eq. (\ref{eq:pt2.3}) we deduce $u^{q} = -\frac{C}{u^2+Au}$ and $u^{q^2} =
\frac{(u^2+Au)^2}{A(u^2+Au)-C}$.
If $u^2+Au=0$, then since $u\neq0$ and substituting $u=-A$ into Eq. (\ref{eq:pt2.3}) leads to $C=0$, which is a contradiction.
If $A(u^2+Au)-C=0$, then we obtain $u=-A$, which is impossible.

Note that $u^{q^2+q+1}=1$ we have
\begin{eqnarray}\label{eq:pt2.6}
u(u+A)(Cu^3 + A(C+1)u^2 + A^{2}u - C)= 0.
\end{eqnarray}
Since $u\neq0$, $u+A\neq0$ and $m_{A,C}(u)=Cu^3 + A(C+1)u^2 + A^{2}u - C$ has no roots in $\mu_{q^2+q+1}$, hence Eq. (\ref{eq:pt2.6}) has no solution.

If $a\neq0$, we show that Eq. (\ref{eq:pt2.1}) has one nonzero solution. Let $y=x^{q}$, $z=y^{q}$, $b=a^{q}$ and $c=b^{q}$, then we obtain the system of equations
\begin{numcases}{}
yz + Azx + Cx^2 -ax =0,  \label{eq:pt2.11}\\
zx + Axy + Cy^2 -by =0,  \label{eq:pt2.11'}\\
xy + Ayz + Cz^2 -cz =0.  \label{eq:pt2.11''}
\end{numcases}
Eliminating $z$ by (\ref{eq:pt2.11}) and (\ref{eq:pt2.11'}), we have
\begin{eqnarray*}
f_{1} \triangleq \xi_{1}y^3 + \xi_{2}y^2 + \xi_{3}y + \xi_{4} =0,
\end{eqnarray*}
where $\xi_{i}:= \xi_{i}(x, A, C, a, b, c)$, $i= 1, 2, 3, 4$.

Furthermore, by Eqs. (\ref{eq:pt2.11}) and (\ref{eq:pt2.11''}) we get
\begin{eqnarray*}
f_{2} \triangleq \eta_{1}y^3 + \eta_{2}y^2 + \eta_{3}y + \eta_{4}=0,
\end{eqnarray*}
where $\eta_{i}:= \eta_{i}(x, A, C, a, b, c)$, $i= 1, 2, 3, 4$.

Computing the resultant of $f_{1}$ and $f_{2}$ with respect for $y$ and recalling that $C^{2}= C-1$, we deduce
\begin{eqnarray*}
x^{2}(Cx - a)^{3}(\alpha x +\beta)=0,
\end{eqnarray*}
where $\alpha := \alpha(A, C, a, b, c)$, $\beta := \beta(A, C, a, b, c)$. Since $a\neq0$ which implies $c\neq0$ , thus $x=0$ is not a solution of the above equation. If $Cx - a=0$, then plugging $x=\frac{a}{C}$ into Eq. (\ref{eq:pt2.1}) derives $-A^3=1$, which is a contradiction. Hence the above equation has at most a solution $x=-\frac{\beta}{\alpha}$.
Therefore Eq. (\ref{eq:pt2.1}) has at most a solution and we complete the proof.
\qed

\begin{theorem}\label{th:pt3}
Let $B, C\in\mathbb{F}_{q}$ with $C^2-C+1=0$, $B^{3}\neq-1$ and $m_{B,C}(t)=Ct^3 -B^2t^2 - BCt - Bt - C\in \mathbb{F}_{q}[t]$ has no roots in $\mu_{q^2+q+1}$. Then
\begin{displaymath}
f(x)= x^{q^2+q-1} + Bx^{q} + Cx
\end{displaymath}
is a permutation trinomial over $\mathbb{F}_{q^3}$.
\end{theorem}

\pf We show that for each $a \in\mathbb{F}_{q^3}$, the equation
\begin{eqnarray}\label{eq:pt3.1}
x^{q^2+q-1} + Bx^{q} + Cx = a
\end{eqnarray}
has at most a solution in $\mathbb{F}_{q^3}$. When $a=0$, we need to prove that $f(x)=0$ has the unique solution $x=0$.  Suppose that $x\neq0$ is one solution of equation
\begin{eqnarray}\label{eq:pt3.2}
x^{q^2+q-1} + Bx^{q} + Cx =0.
\end{eqnarray}
Setting $u=x^{q-1}$ implies $u^{q^2+q+1}=1$. In the following we verify
\begin{eqnarray}\label{eq:pt3.3}
u^{q+2} + Bu + C =0
\end{eqnarray}
has no solution.
From Eq. (\ref{eq:pt3.3}) we deduce $u^{q} = -\frac{Bu+C}{u^2}$ and $u^{q^2} = \frac{B(Bu+C)u^2-Cu^4}{(Bu+C)^2}$.
If $Bu+C=0$, then substituting $u=-\frac{C}{B}$ into Eq. (\ref{eq:pt3.3}) leads to $C=0$, which is a contradiction.

Note that $u^{q^2+q+1}=1$ we have
\begin{eqnarray}\label{eq:pt3.6}
u^2(Bu+C)(Cu^3 -B^2u^2 - BCu - Bu - C)= 0,
\end{eqnarray}
since $u\neq0$, $Bu+C\neq0$ and $m_{B,C}(u)=Cu^3 -B^2u^2 - BCu - Bu - C$ has no roots in $\mu_{q^2+q+1}$, thus Eq. (\ref{eq:pt3.6}) has no solution.

If $a\neq0$, we show that Eq. (\ref{eq:pt3.1}) has one nonzero solution. Let $y=x^{q}$, $z=y^{q}$, $b=a^{q}$ and $c=b^{q}$, then we obtain the system of equations
\begin{numcases}{}
yz + Byx + Cx^2 -ax =0,  \label{eq:pt3.11}\\
zx + Bzy + Cy^2 -by =0,  \label{eq:pt3.11'}\\
xy + Bxz + Cz^2 -cz =0.  \label{eq:pt3.11''}
\end{numcases}
Eliminating $z$ by (\ref{eq:pt3.11}) and (\ref{eq:pt3.11'}), we have
\begin{eqnarray*}
f_{1} \triangleq \xi_{1}y^3 + \xi_{2}y^2 + \xi_{3}y + \xi_{4} =0,
\end{eqnarray*}
where $\xi_{i}:= \xi_{i}(x, B, C, a, b, c)$, $i= 1, 2, 3, 4$.

Furthermore, by Eqs. (\ref{eq:pt3.11}) and (\ref{eq:pt3.11''}) we get
\begin{eqnarray*}
f_{2} \triangleq \eta_{1}y^3 + \eta_{2}y^2 + \eta_{3}y + \eta_{4}=0,
\end{eqnarray*}
where $\eta_{i}:= \eta_{i}(x, B, C, a, b, c)$, $i= 1, 2, 3, 4$.

Computing the resultant of $f_{1}$ and $f_{2}$ with respect for $y$ and recalling that $C^{2}= C-1$, we deduce
\begin{eqnarray*}
x^{2}(Cx - a)^{3}(\alpha x +\beta)=0,
\end{eqnarray*}
where $\alpha := \alpha(B, C, a, b, c)$, $\beta := \beta(B, C, a, b, c)$. Since $a\neq0$, thus $x=0$ is not a solution of the above equation. If $Cx - a=0$, then plugging $x=\frac{a}{C}$ into Eq. (\ref{eq:pt3.1}) derives $-B^3=1$, which is a contradiction. Therefore the above equation has at most a solution $x=-\frac{\beta}{\alpha}$ and
thereby Eq. (\ref{eq:pt3.1}) has at most a solution and the proof is completed.
\qed

\subsection{Four classes of permutation quadrinomials }\label{pd}
In the subsection, we propose four classes of permutation quadrinomials of the form $f(x)= x^{d_{1}} + Ax^{d_{2}} + Bx^{d_{3}} + Cx\in \mathbb{F}_{q}[x]$.

\begin{theorem}\label{th:pq1}
Let $A, B, C \in\mathbb{F}_{q}$ with $A, B\neq0$, $A^{3}\neq-1$, $B^{3}\neq-1$, $C^{2}-AB-C+1=0$, $AB-C\neq 0$, $A + B + 1\neq0$, $A^2 - AB - A + B^2 - B + 1\neq0$ and $m_{A,B,C}(t)=Ct^3 + (AC + A - B^2)t^2 + (A^2 - BC - B)t - C\in \mathbb{F}_{q}[t]$ has no roots in $\mu_{q^2+q+1}$. Then
\begin{displaymath}\label{eq:cpq1.1}
f(x)= x^{q^2+q-1} + Ax^{q^2} + Bx^{q} + Cx
\end{displaymath}
is a permutation quadrinomial over $\mathbb{F}_{q^3}$.
\end{theorem}

\pf We prove that for each $a \in\mathbb{F}_{q^3}$, the equation
\begin{eqnarray}\label{eq:q1.1}
x^{q^2+q-1} + Ax^{q^2} + Bx^{q} + Cx = a
\end{eqnarray}
has at most a solution in $\mathbb{F}_{q^3}$.  We first show $f(x)=0$ only has a solution $x=0$. Assume that $x\neq0$ is one solution of equation
\begin{eqnarray}\label{eq:q1.2}
x^{q^2+q-1} + Ax^{q^2} + Bx^{q} + Cx =0.
\end{eqnarray}
Setting $u=x^{q-1}$ implies $u^{q^2+q+1}=1$. It is necessary to prove that
\begin{eqnarray}\label{eq:q1.3}
u^{q+2} + Au^{q+1} + Bu + C =0
\end{eqnarray}
has no solution.
From Eq. (\ref{eq:q1.3}) we have $u^{q} = -\frac{Bu + C}{u^2+Au}$ and $$u^{q^2} = \frac{-Cu^4+(B^{2}-2AC)u^3 + (BC+AB^{2}-A^{2}C)u^2 +ABCu }{(Bu+C)^{2} - A(Bu+C)(u^2+Au)}.$$
Note that $u^2+Au=0$ implies $u=-A$, substituting $u=-A$ into Eq. (\ref{eq:q1.3}) obtains $AB-C=0$, which is a contradiction.
Furthermore, if $(Bu+C)^{2} - A(Bu+C)(u^2+Au)=0$, then we obtain $u=-\frac{C}{B}$ or $Au^{2}+(A^{2}-B)u-C=0$. Substituting $u=-\frac{C}{B}$ into Eq. (\ref{eq:q1.3}) results in $AB-C=0$, which is a contradiction. Combining $Au^{2}+(A^{2}-B)u-C=0$ and Eq. (\ref{eq:q1.3}) leads to $(Bu+C)(u^q+A)=0$, which means that $Bu+C=0$ or $u^q+A=0$, from the two equations we also obtain $AB-C=0$, which is not possible.

Since $u^{q^2+q+1}=1$ and $AB=C^{2}-C+1$, we deduce
\begin{eqnarray}\label{eq:q1.4}
u(u + A)(Bu + C)(Cu^3 + (AC + A - B^2)u^2 + (A^2 - BC - B)u - C)=0,
\end{eqnarray}
from the above discussion we obtain $u=0$, $u =-A$ and $u =-\frac{C}{B}$ are not solutions of Eq. (\ref{eq:q1.4}), therefore we have
\begin{eqnarray*}
m_{A, B, C}(u)= Cu^3 + (AC + A - B^2)u^2 + (A^2 - BC - B)u - C =0,
\end{eqnarray*}
which contradicts $m_{A, B, C}(u)=0 $ has no solution in $\mu_{q^2+q+1}$ and hence $f(x)=0$ has a unique solution $x=0$.

If $a\neq0$, we show that Eq. (\ref{eq:q1.1}) has one nonzero solution. Let $y=x^{q}$, $z=y^{q}$, $b=a^{q}$ and $c=b^{q}$, then we obtain the system of equations
\begin{numcases}{}
yz + Axz + Bxy+ Cx^2 -ax =0,  \label{eq:q1.11}\\
zx + Ayx + Byz+ Cy^2 -by =0,  \label{eq:q1.11'}\\
xy + Azy + Bzx+ Cz^2 -cz  =0.  \label{eq:q1.11''}
\end{numcases}
Eliminating the indeterminate $z$ by (\ref{eq:q1.11}) and (\ref{eq:q1.11'}), we have
\begin{eqnarray*}
f_{1} \triangleq \xi_{1}y^3 + \xi_{2}y^2 + \xi_{3}y + \xi_{4}=0,
\end{eqnarray*}
where $\xi_{i}:= \xi_{i}(x, A, B, C, a, b, c)$, $i= 1, 2, 3, 4$.

Furthermore, by Eqs. (\ref{eq:q1.11}) and (\ref{eq:q1.11''}) we get
\begin{eqnarray*}
f_{2} \triangleq \eta_{1}y^3 + \eta_{2}y^2 + \eta_{3}y + \eta_{4}=0,
\end{eqnarray*}
where $\eta_{i}:= \eta_{i}(x, A, B, C, a, b, c)$, $i= 1, 2, 3, 4$.

Computing the resultant of $f_{1}$ and $f_{2}$ with respect for $y$ and simplifying the resultant with $AB=C^{2}-C+1$ obtains
\begin{eqnarray*}
C^{2}x^{2}(Cx-a)((C-1)^{2}x+a)^{2}(\alpha x +\beta)=0,
\end{eqnarray*}
where $\alpha := \alpha(A, B, C, a, b, c)$, $\beta := \beta(A, B, C, a, b, c)$. Since $a\neq0$, hence $x=0$ is not a solution of the above equation.

If $Cx-a=0$, then from Eqs. (\ref{eq:q1.11}), (\ref{eq:q1.11'}) and (\ref{eq:q1.11''})  we have $B(A + B + 1)(A^2 - AB - A + B^2 - B + 1)x=0$,  because $B\neq0$, $A + B + 1\neq0$ and $A^2 - AB - A + B^2 - B + 1\neq0$ means $x=0$. It is a contradiction.

If $(C-1)^{2}x+a=0$, then we have $x=-\frac{a}{(C-1)^{2}}$ which implies $y=-\frac{b}{(C-1)^{2}}$ and $z=-\frac{c}{(C-1)^{2}}$. Substituting them into Eq. (\ref{eq:q1.11}) and using $C^{2}=AB+C-1$, we deduce $(Aa+b)(Ba+c)=0$. Suppose that $Ba+c=0$, we derive $c =-Ba$ which implies $a=-Bb$ and $b=-Bc$. Therefore
$c=-Ba=-B(-Bb)=B^2(-Bc)=-B^3c$, yielding $-B^3=1$ which contradicts $B^3+1\neq 0$. The same should happen from $Aa+b=0$.
Hence Eq. (\ref{eq:q1.1}) has at most a solution $x=-\frac{\beta}{\alpha}$ and we complete the proof.
\qed

\begin{theorem}\label{th:pq3}
Let $A, B, C \in\mathbb{F}_{q}$ with $C\neq 0$, $A+B+1\neq0$, $A^3+AB+1\neq0$, $A^2-AB-A+B^2-B+1\neq0$, $A^{3}-ABC+AB+C^{2}-C+1=0$
and $m_{A,B,C}(t)=(AB - C)t^3 + (-A^2C + AB^2 - BC - B)t^2 + (A^2B - AC^2 - AC - B^2)t - AB + C \in \mathbb{F}_{q}[t]$ has no roots in $\mu_{q^2+q+1}$. Then
\begin{displaymath}\label{eq:cpq3.1}
f(x)= x^{q^2+q-1} + Ax^{q^2-q+1} + Bx^{q^2} + Cx
\end{displaymath}
is a permutation quadrinomial over $\mathbb{F}_{q^3}$.
\end{theorem}

\pf We show that for each $a \in\mathbb{F}_{q^3}$, the equation
\begin{eqnarray}\label{eq:q3.1}
x^{q^2+q-1} + Ax^{q^2-q+1} + Bx^{q^2} + Cx = a
\end{eqnarray}
has at most a solution in $\mathbb{F}_{q^3}$.  We prove $x=0$ if and only if $a=0$. Suppose that $x\neq0$ is one solution of equation
\begin{eqnarray}\label{eq:q3.2}
x^{q^2+q-1} + Ax^{q^2-q+1} + Bx^{q^2} + Cx =0.
\end{eqnarray}
Setting $u=x^{q-1}$ implies $u^{q^2+q+1}=1$. It is necessary to prove that
\begin{eqnarray}\label{eq:q3.3}
u^{q+2} + Au^{q} + Bu^{q+1} + C =0
\end{eqnarray}
has no solution.
From Eq. (\ref{eq:q3.3}) we have $u^{q} = -\frac{C}{u^2+Bu + A}$ and
$$u^{q^2} = \frac{C(u^2+Bu+A)^{2} }{-A(u^2+Bu+A)^{2}+BC(u^2+Bu+A)-C^{2}}.$$
If $u^2+Bu+A=0$, then from Eq. (\ref{eq:q3.3}) we have $C=0$, which is impossible.
If $-A(u^2+Bu+A)^{2}+BC(u^2+Bu+A)-C^{2}=0$, then using the equation $A^{3}= ABC-AB- C^{2}+C-1$ and we obtain
\begin{eqnarray}\label{eq:q3d}
Au^4+2ABu^3+(2A^{2}+AB^{2}-BC)u^2 + (2A^{2}B-B^{2}C)u -AB+C-1=0.
\end{eqnarray}
Raising Eq. (\ref{eq:q3d}) to the $q$-th power gets
\begin{eqnarray}\label{eq:q3dq}
&&(AB - C + 1)u^8 + (4AB^2 - 4BC + 4B)u^7 + (2A^2BC + 4A^2B + 6AB^3 - 4AC \nonumber \\
&+& 4A - B^2C^2 - 6B^2C + 6B^2)u^6 + B(6A^2BC + 12A^2B + 4AB^3 - 12AC \nonumber \\
&+& 12A - 3B^2C^2 - 4B^2C + 4B^2)u^5 + (6 A^3 B C + 6 A^3 B + 6 A^2 B^3 C \nonumber \\
&+& 12 A^2 B^3 - 2 A^2 C^2 - 6 A^2 C + 6 A^2 + A B^5 - 4 A B^2 C^2 - 12 A B^2 C \nonumber \\
&+& 12 A B^2 - 3 B^4 C^2 - B^4 C + B^4 + B C^3)u^4 + B(12 A^3 B C + 12 A^3 B \nonumber \\
&+& 2 A^2 B^3 C + 4 A^2 B^3 - 4 A^2 C^2 - 12 A^2 C + 12 A^2 - 8 A B^2 C^2 - 4 A B^2 C \nonumber \\
&+& 4 A B^2 - B^4 C^2 + 2 B C^3)u^3 + (6 A^4 B C + 4 A^4 B + 6 A^3 B^3 C + 6 A^3 B^3 \nonumber \\
&-& 4 A^3 C^2 - 4 A^3 C + 4 A^3 - 7 A^2 B^2 C^2 - 6 A^2 B^2 C + 6 A^2 B^2 - 4 A B^4 C^2 + 4 A B C^3 \nonumber \\
&+& B^3 C^3)u^2 + AB(6 A^3 B C + 4 A^3 B - 4 A^2 C^2 - 4 A^2 C + 4 A^2 - 5 A B^2 C^2 + 4 B C^3)u \nonumber \\
&+& 2 A^5 B C + A^5 B - 2 A^4 C^2 - A^4 C + A^4 - 2 A^3 B^2 C^2 + 3 A^2 B C^3 - A C^4=0.
\end{eqnarray}
Combining Eqs. (\ref{eq:q3d}) and (\ref{eq:q3dq}) and $A^{3}= ABC-AB- C^{2}+C-1$ leads to
\begin{eqnarray*}
(A^3 + AB  + 1)^{24}=0.
\end{eqnarray*}
Since $A^3 + AB  + 1\neq0$ and therefore Eq. (\ref{eq:q3d}) has no solution, which means $-A(u^2+Bu+A)^{2}+BC(u^2+Bu+A)-C^{2}\neq0$.

By $u^{q^2+q+1}=1$ we deduce
\begin{eqnarray*}
&&(u^2+Bu+A)[ Au^4+(2AB-C^{2})u^3 + (2A^2 + AB^{2} - BC^2 - BC)u^2 \nonumber \\
&+&(2A^{2}B - AC^2 -B^{2}C)u  + A^{3}-ABC+C^{2}]=0.
\end{eqnarray*}
Since $C\neq0$,  we have $u^2+Bu+C\neq0$ and $A^{3}-ABC+AB+ C^{2}-C+1=0$,
\begin{eqnarray}\label{eq:q3.4}
&& Au^4+(2AB-C^{2})u^3 + (2A^2 + AB^{2} - BC^2 - BC)u^2 +(2A^{2}B - AC^2 -B^{2}C)u  \nonumber \\
&-& AB + C - 1=0.
\end{eqnarray}
Raising Eq. (\ref{eq:q3.4}) to the $q$-th power results in
\begin{eqnarray}\label{eq:q3.5}
&& (AB - C + 1)u^8 + (4AB^2 - 4BC + 4B)u^7 + (2A^2BC + 4 A^2 B + 6 AB^3 - AC^3 - 4AC \nonumber \\
&+& 4A - B^2C^2 - 6B^2C + 6B^2)u^6 + B(6A^2BC + 12A^2B + 4AB^3 - 3AC^3 - 12AC + 12A \nonumber \\
&-& 3B^2C^2 - 4B^2C + 4B^2)u^5
+ (6A^3BC + 6A^3B + 6A^2B^3C + 12A^2B^3 - 3A^2C^3 - 2A^2C^2 \nonumber \\
&-& 6A^2C + 6A^2 + AB^5 - 3AB^2C^3 - 4A B^2 C^2 - 12 A B^2 C + 12 A B^2 - 3B^4C^2 - B^4C  \nonumber \\
&+& B^4 + BC^4 + BC^3)u^4
+ B(12 A^3 B C + 12 A^3 B + 2 A^2 B^3 C + 4 A^2 B^3 - 6 A^2 C^3 - 4 A^2C^2 \nonumber \\
&-& 12A^2C + 12A^2 - AB^2C^3 - 8AB^2C^2 - 4 AB^2C + 4AB^2 - B^4C^2 + 2BC^4 + 2BC^3)u^3  \nonumber \\
&+& (6 A^4 B C + 4 A^4 B +6 A^3 B^3 C + 6 A^3 B^3 - 3 A^3 C^3 - 4 A^3 C^2 - 4 A^3 C + 4A^3 - 3A^2B^2C^3 \nonumber \\
&-& 7A^2B^2C^2 - 6A^2B^2C + 6A^2B^2 - 4AB^4C^2 + 2 A B C^4 + 4 A B C^3 + B^3C^4 + B^3C^3 - C^5)u^2 \nonumber \\
&+& B(6A^4 B C + 4A^4 B - 3A^3 C^3 -4A^3C^2 -4A^3C + 4 A^3 - 5A^2 B^2 C^2 + 2ABC^4 + 4ABC^3 \nonumber \\
&-& C^5)u + 2A^5BC + A^5B - A^4C^3 - 2A^4C^2 - A^4C + A^4 - 2A^3B^2 C^2 + A^2BC^4 + 3A^2BC^3 \nonumber \\
&-& AC^5 - AC^4=0.
 \end{eqnarray}
Combinging Eqs. (\ref{eq:q3.4}) and (\ref{eq:q3.5}) and $A^{3}= ABC-AB- C^{2}+C-1$, we derive
\begin{eqnarray*}
C^5((AB - C)u^3 + (-A^2C + AB^2 - BC - B)u^2 + (A^2B -
        AC^2 - AC - B^2)u - AB + C)=0,
\end{eqnarray*}
since $C\neq0$, we have $m_{A, B, C}(u)=(AB - C)u^3 + (-A^2C + AB^2 - BC - B)u^2 + (A^2B - AC^2 - AC - B^2)u - AB + C=0$, which contradicts $m_{A, B, C}(u)=0 $ has no solution in $\mu_{q^2+q+1}$ and hence $f(x)=0$ has a unique solution $x=0$.

If $a\neq0$, we show that Eq. (\ref{eq:q3.1}) has one nonzero solution. Let $y=x^{q}$, $z=y^{q}$, $b=a^{q}$ and $c=b^{q}$, then we obtain the system of equations
\begin{numcases}{}
y^2z + Ax^2z + Bxyz+ Cx^2y -axy =0,  \label{eq:q3.11}\\
z^2x + Ay^2x + Bxyz+ Cy^2z -byz =0,  \label{eq:q3.11'}\\
x^2y + Az^2y + Bxyz+ Cz^2x -czx  =0.  \label{eq:q3.11''}
\end{numcases}
Eliminating the indeterminate $z$ by (\ref{eq:q3.11}) and (\ref{eq:q3.11'}), we have
\begin{eqnarray*}
f_{1} \triangleq \xi_{1}y^4 + \xi_{2}y^3 + \xi_{3}y^2 + \xi_{4}y + \xi_{5}=0,
\end{eqnarray*}
where $\xi_{i}:= \xi_{i}(x, A, B, C, a, b, c)$, $i= 1, 2, 3, 4, 5$.

Furthermore, by Eqs. (\ref{eq:q3.11}) and (\ref{eq:q3.11''}) we get
\begin{eqnarray*}
f_{2} \triangleq \eta_{1}y^4 + \eta_{2}y^3 + \eta_{3}y^2 + \eta_{4}y + \eta_{5}=0,
\end{eqnarray*}
where $\eta_{i}:= \eta_{i}(x, A, B, C, a, b, c)$, $i= 1, 2, 3, 4, 5$.

Computing the resultant of $f_{1}$ and $f_{2}$ with respect for $y$ and recalling that $A^{3}= ABC-AB- C^{2}+C-1$, we derive
\begin{eqnarray*}
x^{4}(Cx - a)^{8}(\alpha x +\beta)=0,
\end{eqnarray*}
where $\alpha := \alpha(A, B, C, a, b, c)$, $\beta := \beta(A, B, C, a, b, c)$. Since $a\neq0$, hence $x=0$ is not a solution of the above equation. If $x=\frac{a}{C}$ is a solution, then from Eqs. (\ref{eq:q3.11}), (\ref{eq:q3.11'}) and (\ref{eq:q3.11''}) we obtain $(A+B+1)(A^3+AB+1)^{3}(A^2-AB-A+B^2-B+1)x^8=0$, since $A+B+1\neq0$, $A^3+AB+1\neq0$ and  $A^2-AB-A+B^2-B+1\neq0$, thus $x=0$, which is a contradiction.
Therefore Eq. (\ref{eq:q3.1}) has at most a solution $x=-\frac{\beta}{\alpha}$ and we complete the proof.
\qed

Let $r(A, C) = (3A^9C^2 - 2A^9C + 11A^6C^4 - 17A^6C^3 + 7A^6C^2 + 6A^6C - 4A^6 + A^3C^8 - A^3C^7 + 6A^3C^6 - 17A^3C^5 + 23A^3C^4 - 12A^3C^3 + 3A^3C^2 + C^8 - 3C^7 + 6C^6 - 7C^5 + 6C^4 - 3C^3 + C^2)(A^9C^3 - 4A^9C^2 + 4A^9C + 5A^6C^5 - 23A^6C^4 + 33A^6C^3 - 10A^6C^2 - 12A^6C + 8A^6 - A^3C^8 + 10A^3C^7 - 36A^3C^6 + 57A^3C^5 - 38A^3C^4 - A^3C^3 + 16A^3C^2 - 8A^3C - C^10 + 6C^9 - 17C^8 + 27C^7 - 25C^6 + 10C^5 + 3C^4 - 5C^3 + 2C^2)$,

$r_{1}(A,C)=(22A^{27}C - 14A^{27} - 37A^{24}C^5 + 99A^{24}C^4 + 1210A^{24}C^3 -
        2488A^{24}C^2 + 1600A^{24}C - 336A^{24} + 15A^{21}C^9 - 87A^{21}C^8 -
        1271A^{21}C^7 + 5737A^{21}C^6 + 3800A^{21}C^5 - 35038A^{21}C^4 +
        50748A^{21}C^3 - 32984A^{21}C^2 + 10360A^{21}C - 1288A^{21} +
        12A^{18}C^{12} + 290A^{18}C^{11} - 2516A^{18}C^{10} - 605A^{18}C^9 +
        32787A^{18}C^8 - 69632A^{18}C^7 + 44484A^{18}C^6 + 5966A^{18}C^5 -
        7086A^{18}C^4 - 12012A^{18}C^3 + 10920A^{18}C^2 - 2912A^{18}C + 208A^{18}
        + 156A^{15}C^{14} - 276A^{15}C^{13} - 4186A^{15}C^{12} + 22731A^{15}C^{11} -
        75149A^{15}C^{10} + 171008A^{15}C^9 - 33196A^{15}C^8 - 942781A^{15}C^7 +
        2546731A^{15}C^6 - 3326620A^{15}C^5 + 2545906A^{15}C^4 -
        1186900A^{15}C^3 + 328328A^{15}C^2 - 48620A^{15}C + 2860A^{15} +
        12A^{12}C^{16} - 1432A^{12}C^{15} + 19052A^{12}C^{14} - 88879A^{12}C^{13} +
        88879A^{12}C^{12} + 621676A^{12}C^{11} - 2430430A^{12}C^{10} +
        3274376A^{12}C^9 + 139178A^{12}C^8 - 6582658A^{12}C^7 +
        10211906A^{12}C^6 - 8340288A^{12}C^5 + 4159232A^{12}C^4 -
        1277848A^{12}C^3 + 225368A^{12}C^2 - 18304A^{12}C + 208A^{12} -
        444A^9C^{18} + 5289A^9C^{17} - 16063A^9C^{16} - 51418A^9C^{15} +
        455342A^9C^{14} - 1098227A^9C^{13} + 516193A^9C^{12} + 2843722A^9C^{11}
        - 6952122A^9C^{10} + 7234639A^9C^9 - 3222583A^9C^8 - 513169A^9C^7
        + 1104863A^9C^6 - 90992A^9C^5 - 447786A^9C^4 + 322868A^9C^3 -
        107016A^9C^2 + 18200A^9C - 1288A^9 + 168A^6C^20 + 162A^6C^{19} -
        15656A^6C^{18} + 93907A^6C^{17} - 214051A^6C^{16} + 15156A^6C^{15} +
        1037482A^6C^{14} - 2608867A^6C^{13} + 3126827A^6C^{12} -
        1416714A^6C^{11} - 1893680A^6C^{10} + 5178859A^6C^9 - 7077049A^6C^8
        + 6924516A^6C^7 - 5076720A^6C^6 + 2770090A^6C^5 - 1102178A^6C^4
        + 309708A^6C^3 - 58184A^6C^2 + 6560A^6C - 336A^6 + 96A^3C^{22} -
        1424A^3C^{21} + 7744A^3C^{20} - 15684A^3C^{19} - 21180A^3C^{18} +
        191679A^3C^{17} - 466253A^3C^{16} + 506054A^3C^{15} + 135174A^3C^{14} -
        1419473A^3C^{13} + 2575061A^3C^{12} - 2690813A^3C^{11} +
        1610495A^3C^{10} - 111070A^3C^9 - 862278A^3C^8 + 1007205A^3C^7 -
        672531A^3C^6 + 307316A^3C^5 - 99690A^3C^4 + 22756A^3C^3 -
        3496A^3C^2 + 326A^3C - 14A^3 - 32C^{23} + 480C^{22} - 3024C^{21} +
        10328C^{20} - 19486C^{19} + 12448C^{18} + 35461C^{17} - 119645C^{16} +
        184024C^{15} - 161830C^{14} + 48906C^{13} + 81918C^{12} - 150300C^{11} +
        138410C^{10} - 86294C^9 + 38780C^8 - 12674C^7 + 2954C^6 - 467C^5 +
        45C^4 - 2C^3)$.

\begin{theorem}\label{th:pq4}
Let $A, B, C\in\mathbb{F}_{q}$ with $A\neq0$, $C\neq0, 1$, $A^3 + AB  + 1\neq0$, $B^2=4A$, $A^{3}- ABC+AB+C^{2}-C+1=0$, $r(A,C)\neq0$ and $r_{1}(A,C)\neq0$. Then
\begin{displaymath}\label{eq:cpq4.1}
f(x)= x^{q^2+q-1} + Ax^{q^3-q^2+q} + Bx^{q} + Cx
\end{displaymath}
is a permutation quadrinomial over $\mathbb{F}_{q^3}$.
\end{theorem}

\pf We prove that for each $a \in\mathbb{F}_{q^3}$, the equation
\begin{eqnarray}\label{eq:q4.1}
x^{q^2+q-1} + Ax^{q^3-q^2+q} + Bx^{q} + Cx = a
\end{eqnarray}
has at most a solution in $\mathbb{F}_{q^3}$.  We demonstrate $f(x)=0$ only has a solution $x=0$. Suppose that  $x\neq0$ is one solution of equation
\begin{eqnarray}\label{eq:q4.2}
x^{q^2+q-1} + Ax^{q^3-q^2+q} + Bx^{q} + Cx =0.
\end{eqnarray}
Setting $u=x^{q-1}$ implies $u^{q^2+q+1}=1$. It is need to prove that
\begin{eqnarray}\label{eq:q4.3}
u^{q+2} + Au^{q^2+1} + Bu + C =0
\end{eqnarray}
has no solution.
Raising both sides of Eq. (\ref{eq:q4.3}) to the $q$-th power and then multiplying by $u$ leads to
\begin{eqnarray}\label{eq:q4.4}
u^{q} + Au^{q+2} + Bu^{q+1} + Cu =0,
\end{eqnarray}
from Eq. (\ref{eq:q4.3}) we have $u^{q} = -\frac{Cu}{Au^2+Bu + 1}$ and $$u^{q^2} = \frac{C^2u(Au^{2}+Bu +1) }{(Au^{2}+Bu +1)^{2} - BCu(Au^{2}+Bu +1)+AC^{2}u^{2}}.$$
If $Au^2+Bu + 1=0$ implies $C=0$ from Eq. (\ref{eq:q4.4}), which is a contradiction since $C\neq0$.
Furthermore, if $(Au^2+Bu + 1)^{2} - BCu(Au^{2}+Bu +1)+AC^{2}u^{2}=0$, then we obtain
\begin{eqnarray}\label{eq:q4d}
A^{2}u^4+(2AB-ABC)u^3+(AC^{2}+2A-B^{2}C+B^{2})u^2 + (2B-BC)u +1=0.
\end{eqnarray}
Raising Eq. (\ref{eq:q4d}) to the $q$-th power with $A^{3}= ABC-AB- C^{2}+C-1$ gets
\begin{eqnarray}\label{eq:q4dq}
&& A^4u^8 + (A^3BC^2 - 2A^3BC + 4A^3B)u^7 \nonumber \\
&+& A^2(AC^4 + 2AC^2 + 4A - B^2C^3 + 4B^2C^2 - 6B^2C + 6B^2)u^6 \nonumber \\
&+& AB(3AC^4 - 2AC^3 + 7AC^2 - 6AC + 12A - 2B^2C^3 + 5B^2C^2 - 6B^2C + 4B^2)u^5 \nonumber \\
&+& (3 A^2 C^4 + 4 A^2 C^2 + 6 A^2 + 2 A B^2 C^4 - 4 A B^2 C^3 + 10 A B^2 C^2 - 12 A B^2 C \nonumber \\
&+& 12 A B^2 - B^4 C^3 + 2 B^4 C^2 - 2 B^4 C + B^4)u^4 \nonumber \\
&+& B(3AC^4 - 2AC^3 + 7AC^2 - 6AC + 12A - 2B^2C^3 + 5B^2C^2 - 6B^2C + 4B^2)u^3 \nonumber \\
&+& (AC^4 + 2AC^2 + 4A - B^2C^3 + 4B^2C^2 - 6B^2C + 6B^2)u^2 \nonumber \\
&+& (BC^2 - 2BC + 4B)u + 1=0.
\end{eqnarray}
Combining Eqs. (\ref{eq:q4d}) and (\ref{eq:q4dq}) eliminates $u$ we have
\begin{eqnarray*}
A^{32}(A^3 + AB  + 1)^{24}=0,
\end{eqnarray*}
since $A \neq0$ and $A^3 + AB  + 1\neq0$ and therefore  Eq. (\ref{eq:q4d}) has no solution, which means $(Au^2+Bu + 1)^{2} - BCu(Au^{2}+Bu +1)+AC^{2}u^{2}\neq0$.

Note that $u^{q^2+q+1}=1$ and $A^{3}= ABC-AB- C^{2}+C-1$, we derive
$$(Au^2+Bu+1)[ A^2u^4+(2AB-ABC+C^{3})u^3 + (2A - AC^{2} - B^2C + B^{2})u^2 + (2B - BC)u + 1 ]=0. $$
Since $C\neq0$, hence $u^2+Bu+A\neq0$ and
\begin{eqnarray}\label{eq:q4.5}
A^2u^4+(2AB-ABC+C^{3})u^3 + (2A - AC^{2} - B^2C + B^{2})u^2 + (2B - BC)u + 1 =0.
\end{eqnarray}
Raising Eq. (\ref{eq:q4.5}) to the $q$-th power gets
\begin{eqnarray}\label{eq:q4.5'}
&& A^4u^8 + (A^3BC^2 - 2A^3BC + 4A^3B)u^7 \nonumber \\
&+& A^2(-AC^4 + 2AC^2 + 4A - B^2C^3 + 4B^2C^2 - 6B^2C + 6B^2)u^6  \nonumber \\
&+& A(-ABC^4 - 2ABC^3 + 7ABC^2 - 6ABC + 12AB - 2B^3C^3 + 5B^3C^2 - 6B^3C + 4B^3 - C^6)u^5 \nonumber \\
&+& (-A^2C^4 + 4A^2C^2 + 6A^2 - 4AB^2C^3 + 10AB^2C^2 - 12AB^2C \nonumber \\
&+& 12AB^2 - B^4C^3 + 2B^4C^2 - 2B^4C + B^4 - BC^6)u^4 \nonumber \\
&+& (-3ABC^4 - 2ABC^3 + 7ABC^2 - 6ABC + 12AB - 2B^3C^3 + 5B^3C^2 - 6B^3C + 4B^3 - C^6)u^3 \nonumber \\
&+& (-AC^4 + 2AC^2 + 4A - B^2C^3 + 4B^2C^2 - 6B^2C + 6B^2)u^2 \nonumber \\
&+& (BC^2 - 2BC + 4B)u + 1=0.
\end{eqnarray}
Combining Eqs. (\ref{eq:q4.5}) and (\ref{eq:q4.5'}) results in
\begin{eqnarray*}
2^3A^{16}C^{16}(C - 1)^9\cdot r(A, C)=0,
\end{eqnarray*}
where $r(A, C) = (3A^9C^2 - 2A^9C + 11A^6C^4 - 17A^6C^3 + 7A^6C^2 + 6A^6C - 4A^6 + A^3C^8 - A^3C^7 + 6A^3C^6 - 17A^3C^5 + 23A^3C^4 - 12A^3C^3 + 3A^3C^2 + C^8 - 3C^7 + 6C^6 - 7C^5 + 6C^4 - 3C^3 + C^2)(A^9C^3 - 4A^9C^2 + 4A^9C + 5A^6C^5 - 23A^6C^4 + 33A^6C^3 - 10A^6C^2 - 12A^6C + 8A^6 - A^3C^8 + 10A^3C^7 - 36A^3C^6 + 57A^3C^5 - 38A^3C^4 - A^3C^3 + 16A^3C^2 - 8A^3C - C^10 + 6C^9 - 17C^8 + 27C^7 - 25C^6 + 10C^5 + 3C^4 - 5C^3 + 2C^2)$,
however, $A\neq0$, $C\neq0, 1$ and $r(A, C)\neq0$ and therefore $f(x)=0$ has a unique solution $x=0$.

If $a\neq0$, we show that Eq. (\ref{eq:q4.1}) has one nonzero solution. Let $y=x^{q}$, $z=y^{q}$, $b=a^{q}$ and $c=b^{q}$, then we obtain the system of equations
\begin{numcases}{}
yz^2 + Ax^2y + Bxyz+ Cx^2z -axz =0,  \label{eq:q4.11}\\
zx^2 + Ay^2z + Bxyz+ Cy^2x -byx =0,  \label{eq:q4.11'}\\
xy^2 + Az^2x + Bxyz+ Cz^2y -czy  =0.  \label{eq:q4.11''}
\end{numcases}
Eliminating the indeterminate $z$ by (\ref{eq:q4.11}) and (\ref{eq:q4.11'}), we have
\begin{eqnarray*}
f_{1} \triangleq \xi_{1}y^4 + \xi_{2}y^3 + \xi_{3}y^2 + \xi_{4}y + \xi_{5}=0,
\end{eqnarray*}
where $\xi_{i}:= \xi_{i}(x, A, B, C, a, b, c)$, $i= 1, 2, 3, 4, 5$.

Furthermore, by Eqs. (\ref{eq:q4.11'}) and (\ref{eq:q4.11''}) we get
\begin{eqnarray*}
f_{2} \triangleq \eta_{1}y^4 + \eta_{2}y^3 + \eta_{3}y^2 + \eta_{4}y + \eta_{5}=0,
\end{eqnarray*}
where $\eta_{i}:= \eta_{i}(x, A, B, C, a, b, c)$, $i= 1, 2, 3, 4, 5$.

We compute the resultant of $f_{1}$ and $f_{2}$ with respect for $y$ and recall that $A^{3}= ABC-AB- C^{2}+C-1$,
\begin{eqnarray*}
x^{4}(C^{2}x^2 + bBCx + b^2A)^{4}(\alpha x +\beta)=0,
\end{eqnarray*}
where $\alpha := \alpha(A, B, C, a, b, c)$, $\beta := \beta(A, B, C, a, b, c)$. Since $a\neq0$, hence $x=0$ is not a solution of the above equation.

If $C^{2}x^2 + bBCx + b^2A=0$, then by $B^2=4A$ we have $x=-\frac{bB}{2C}$. Substituting $x=-\frac{bB}{2C}$ into Eqs.  (\ref{eq:q4.11}), (\ref{eq:q4.11'}), (\ref{eq:q4.11''}) and eliminating $b$, $c$ with $A^{3}= ABC-AB- C^{2}+C-1$, we deduce
$$ 2^{19}A^{8}r_{1}(A,C)\cdot a^{16}=0, $$
where $r_{1}(A,C)=(22A^{27}C - 14A^{27} - 37A^{24}C^5 + 99A^{24}C^4 + 1210A^{24}C^3 -
        2488A^{24}C^2 + 1600A^{24}C - 336A^{24} + 15A^{21}C^9 - 87A^{21}C^8 -
        1271A^{21}C^7 + 5737A^{21}C^6 + 3800A^{21}C^5 - 35038A^{21}C^4 +
        50748A^{21}C^3 - 32984A^{21}C^2 + 10360A^{21}C - 1288A^{21} +
        12A^{18}C^{12} + 290A^{18}C^{11} - 2516A^{18}C^{10} - 605A^{18}C^9 +
        32787A^{18}C^8 - 69632A^{18}C^7 + 44484A^{18}C^6 + 5966A^{18}C^5 -
        7086A^{18}C^4 - 12012A^{18}C^3 + 10920A^{18}C^2 - 2912A^{18}C + 208A^{18}
        + 156A^{15}C^{14} - 276A^{15}C^{13} - 4186A^{15}C^{12} + 22731A^{15}C^{11} -
        75149A^{15}C^{10} + 171008A^{15}C^9 - 33196A^{15}C^8 - 942781A^{15}C^7 +
        2546731A^{15}C^6 - 3326620A^{15}C^5 + 2545906A^{15}C^4 -
        1186900A^{15}C^3 + 328328A^{15}C^2 - 48620A^{15}C + 2860A^{15} +
        12A^{12}C^{16} - 1432A^{12}C^{15} + 19052A^{12}C^{14} - 88879A^{12}C^{13} +
        88879A^{12}C^{12} + 621676A^{12}C^{11} - 2430430A^{12}C^{10} +
        3274376A^{12}C^9 + 139178A^{12}C^8 - 6582658A^{12}C^7 +
        10211906A^{12}C^6 - 8340288A^{12}C^5 + 4159232A^{12}C^4 -
        1277848A^{12}C^3 + 225368A^{12}C^2 - 18304A^{12}C + 208A^{12} -
        444A^9C^{18} + 5289A^9C^{17} - 16063A^9C^{16} - 51418A^9C^{15} +
        455342A^9C^{14} - 1098227A^9C^{13} + 516193A^9C^{12} + 2843722A^9C^{11}
        - 6952122A^9C^{10} + 7234639A^9C^9 - 3222583A^9C^8 - 513169A^9C^7
        + 1104863A^9C^6 - 90992A^9C^5 - 447786A^9C^4 + 322868A^9C^3 -
        107016A^9C^2 + 18200A^9C - 1288A^9 + 168A^6C^20 + 162A^6C^{19} -
        15656A^6C^{18} + 93907A^6C^{17} - 214051A^6C^{16} + 15156A^6C^{15} +
        1037482A^6C^{14} - 2608867A^6C^{13} + 3126827A^6C^{12} -
        1416714A^6C^{11} - 1893680A^6C^{10} + 5178859A^6C^9 - 7077049A^6C^8
        + 6924516A^6C^7 - 5076720A^6C^6 + 2770090A^6C^5 - 1102178A^6C^4
        + 309708A^6C^3 - 58184A^6C^2 + 6560A^6C - 336A^6 + 96A^3C^{22} -
        1424A^3C^{21} + 7744A^3C^{20} - 15684A^3C^{19} - 21180A^3C^{18} +
        191679A^3C^{17} - 466253A^3C^{16} + 506054A^3C^{15} + 135174A^3C^{14} -
        1419473A^3C^{13} + 2575061A^3C^{12} - 2690813A^3C^{11} +
        1610495A^3C^{10} - 111070A^3C^9 - 862278A^3C^8 + 1007205A^3C^7 -
        672531A^3C^6 + 307316A^3C^5 - 99690A^3C^4 + 22756A^3C^3 -
        3496A^3C^2 + 326A^3C - 14A^3 - 32C^{23} + 480C^{22} - 3024C^{21} +
        10328C^{20} - 19486C^{19} + 12448C^{18} + 35461C^{17} - 119645C^{16} +
        184024C^{15} - 161830C^{14} + 48906C^{13} + 81918C^{12} - 150300C^{11} +
        138410C^{10} - 86294C^9 + 38780C^8 - 12674C^7 + 2954C^6 - 467C^5 +
        45C^4 - 2C^3)$.

Since $ A\neq0$ and $r_{1}(A,C)\neq0$, we obtain $a=0$ which contradicts $a\neq0$.
Therefore Eq. (\ref{eq:q4.1}) has at most a solution and the proof is completed.
\qed

Let $r_{1}(v,v_{1})=(v + v_{1} + 1)(v^2 + v + 1)^{6}(v^2 - vv_{1} - v + v_{1}^2 - v_{1} + 1)(v^3 + vv_{1} + 1)^{5}(v^4 - 2v^2v_{1} + vv_{1}^3 - v + v_{1}^2)(v^4 - v^2v_{1} + v - v_{1}^2)^{2}(v^4 + 2v^2v_{1} - v + v_{1}^2)$,

$r_{2}(v,v_{1})=(v - v_{1} - 1)(v^2 - v + 1)^{6}(v^2 + vv_{1} + v + v_{1}^2 - v_{1} + 1)(v^3 + vv_{1} - 1)^{5}(v^4 - 2v^2v_{1} - vv_{1}^3 + v + v_{1}^2)(v^4 - v^2v_{1} - v - v_{1}^2)^{2}(v^4 + 2v^2v_{1} + v + v_{1}^2)$.

\begin{theorem}\label{th:pq2}
Let $A, B, C\in\mathbb{F}_{q}$ with $A\neq-1, 0, 1$, $A^2 - A + 1\neq0$, $AB^{2}+C^{2}\neq 0$,
$A^3 - 2A^2B + AB^2 + 1\neq0$ and $A^{3}-A^{2}B+AB^{2}+C^{2}-C+1=0$, $r_{1}(A,B)\neq0$, $r_{2}(A,B)\neq0$, $m_{A, B, C}(t)=\lambda_{1}t^{3} + \lambda_{2}t^2 + \lambda_{3}t + \lambda_{4}\in \mathbb{F}_{q}[t]$ has no roots in $\mu_{q^2+q+1}$, $\lambda_{i}=\lambda_{i}(A, B, C), i=1, 2, 3, 4$, (see Appendix \ref{App}). Then
\begin{displaymath}\label{eq:cpq2.1}
f(x)= x^{q^2+q-1} + Ax^{q^2-q+1} + Bx^{q} + Cx
\end{displaymath}
is a permutation quadrinomial over $\mathbb{F}_{q^3}$.
\end{theorem}

\pf We prove that for each $a \in\mathbb{F}_{q^3}$, the equation
\begin{eqnarray}\label{eq:q2.1}
x^{q^2+q-1} + Ax^{q^2-q+1} + Bx^{q} + Cx = a
\end{eqnarray}
has at most a solution in $\mathbb{F}_{q^3}$.  We demonstrate $x=0$ if and only if $a=0$. Assume that $x\neq0$ is one solution of equation
\begin{eqnarray}\label{eq:q2.2}
x^{q^2+q-1} + Ax^{q^2-q+1} + Bx^{q} + Cx =0.
\end{eqnarray}
Setting $u=x^{q-1}$ implies $u^{q^2+q+1}=1$. It is necessary to prove that
\begin{eqnarray}\label{eq:q2.3}
u^{q+2} + Au^{q} + Bu + C =0
\end{eqnarray}
has no solution.
From Eq. (\ref{eq:q2.3}) we have $u^{q} = -\frac{Bu + C}{u^2+A}$ and $u^{q^2} = \frac{(B^{2}u+BC)(u^2+A) - C(u^2+A)^{2} }{(Bu+C)^{2} + A(u^2+A)^{2}}$.
Note that $u^2+A=0$, substituting it into Eq. (\ref{eq:q2.3}) obtains $AB^{2}+C^{2}=0$, which is a contradiction.
If $(Bu+C)^{2} +A(u^2+A)^{2}=0$, then we obtain
\begin{eqnarray}\label{eq:q2d}
Au^4+(2A^{2}+B^{2})u^2 + 2BCu+A^{3}+C^{2}=0.
\end{eqnarray}
Raising Eq. (\ref{eq:q2d}) to the $q$-th power gets
\begin{eqnarray}\label{eq:q2dq}
&&(A^3 + C^2)u^8 - 2B^2Cu^7 + (4A^4 + 2A^2B^2 + 4AC^2 +B^4 - 2BC^2)u^6  \nonumber \\
&+& (4A^2BC - 6AB^2C + 2B^3C)u^5 + (6A^5 + 4A^3B^2 + 8A^2C^2 + 3AB^4 - 6ABC^2 + B^2C^2 )u^4  \nonumber \\
&+& (8A^3BC - 6A^2B^2C + 8AB^3C)u^3 + (4A^6 + 2A^4B^2 + 8A^3C^2 + A^2B^4 - 6A^2BC^2 + 8AB^2C^2)u^2 \nonumber \\
&+& (4A^4BC - 2A^3B^2C + 2A^2B^3C + 4ABC^3)u \nonumber \\
&+& A^7 + 3A^4C^2 - 2A^3BC^2 + A^2B^2C^2 + AC^4=0.
\end{eqnarray}
By Eqs. (\ref{eq:q2d}) and (\ref{eq:q2dq}) eliminate $u$ we have
\begin{eqnarray*}
(A + 1)^{12}(A^2 - A + 1)^{12}(A^3 - 2A^2B + AB^2 + 1)^{12}=0,
\end{eqnarray*}
since $A\neq-1$, $A^2 - A + 1\neq0$ and $A^3 - 2A^2B + AB^2 + 1\neq0$, which means  Eq. (\ref{eq:q2d}) has no solution, thereby $(Bu+C)^{2} - A(u^2+A)^{2}\neq0$.

Note that $u^{q^2+q+1}=1$, we obtain
\begin{eqnarray*}
&& (u^2+A)[ (A-BC)u^4+(B^{3}-C^{2})u^3 + (2A^2 - ABC + 2B^2C + B^2)u^2 \nonumber \\
&+& (-AC^2 + BC^2 + 2BC)u + A^{3}+C^{2}]=0.
\end{eqnarray*}
Since $AB^{2}+C^{2}\neq0$, we get $u^2+A\neq0$ and by $A^{3}= A^{2}B-AB^{2}- C^{2}+C-1$ we derive
\begin{eqnarray*}
&&(A-BC)u^4+(B^{3}-C^{2})u^3 + (2A^2 - ABC + 2B^2C + B^2)u^2 +(-AC^2 + BC^2 + 2BC)u \nonumber \\
&+& A^{2}B - AB^{2} + C-1=0.
\end{eqnarray*}
Raising both sides of the above equation to the $q$-th power results in
\begin{eqnarray}\label{eq:q2.5}
&& (A^2B - AB^2 + C - 1)u^8 + (ABC^2 - B^2C^2 - 2B^2C)u^7 + (4A^3B - 2A^2B^2 - AB^3C + AC^3 \nonumber \\
&+& 4AC - 4A + 2B^4C + B^4 - BC^3 - 2BC^2)u^6+B(3A^2C^2 + 4A^2C - 5ABC^2 - 6ABC - B^5 \nonumber \\
&+& 5B^2C^2 + 2B^2C)u^5+ (6A^4B - 2A^3B^2 - 2A^2B^3C + 3A^2C^3 + 2A^2C^2 + 6A^2C - 6A^2 \nonumber \\
&+& 4AB^4C + 3AB^4 - 4ABC^3 - 6ABC^2 - 4B^5C + 5B^2C^3 + B^2C^2)u^4 + B(3A^3C^2 + 8A^3C \nonumber \\
&-& 7A^2BC^2 - 6A^2BC - AB^5 + 9AB^2C^2 + 8AB^2C - 7B^3C^2 + 3C^4)u^3+ (4A^5B - 2A^4B^2 \nonumber \\
&-& A^3B^3C + 3A^3C^3 + 4A^3C^2 + 4A^3C - 4A^3 + 2A^2B^4C + A^2B^4 - 5A^2BC^3 - 6A^2BC^2 \nonumber \\
&-& 3AB^5C + 7AB^2C^3 + 8AB^2C^2 - 7B^3C^3 + C^5)u^2+ BC(A^4C + 4A^4 - 3A^3BC - 2A^3B \nonumber \\
&+& 4A^2B^2C + 2A^2B^2 - 3AB^3C + 3AC^3 + 4AC^2 - 4BC^3)u+ A^6B -A^5B^2 + A^4C^3 \nonumber \\
&+& 2A^4C^2 + A^4C - A^4 - 2A^3BC^3 - 2A^3BC^2+ 2A^2B^2C^3 + A^2B^2C^2 - AB^3C^3 \nonumber \\
&+& AC^5 + AC^4 - BC^5 = 0.
\end{eqnarray}
Through a series of computations we obtain $m_{A, B, C}(u)=\lambda_{1}u^{3} + \lambda_{2}u^2 + \lambda_{3}u + \lambda_{4}=0 $, $\lambda_{i}=\lambda_{i}(A, B, C), i=1, 2, 3, 4$, (see Appendix \ref{App}), which contradicts $m_{A, B, C}(u)=0 $ has no solution in $\mu_{q^2+q+1}$ and therefore $f(x)=0$ has a unique solution $x=0$.

If $a\neq0$, we demonstrate that Eq. (\ref{eq:q2.1}) has one nonzero solution. Let $y=x^{q}$, $z=y^{q}$, $b=a^{q}$ and $c=b^{q}$, then we obtain the system of equations
\begin{numcases}{}
y^2z + Ax^2z + Bxy^2+ Cx^2y -axy =0,  \label{eq:q2.11}\\
z^2x + Ay^2x + Byz^2+ Cy^2z -byz =0,  \label{eq:q2.11'}\\
x^2y + Az^2y + Bzx^2+ Cz^2x -czx  =0.  \label{eq:q2.11''}
\end{numcases}
Eliminating the indeterminate $z$ by (\ref{eq:q2.11}) and (\ref{eq:q2.11'}), we have
\begin{eqnarray*}
f_{1} \triangleq \xi_{1}y^4 + \xi_{2}y^3 + \xi_{3}y^2 + \xi_{4}y + \xi_{5}=0,
\end{eqnarray*}
where $\xi_{i}:= \xi_{i}(x, A, B, C, a, b, c)$, $i= 1, 2, 3, 4, 5$.

Furthermore, by Eqs. (\ref{eq:q2.11}) and (\ref{eq:q2.11''}) we get
\begin{eqnarray*}
f_{2} \triangleq \eta_{1}y^4 + \eta_{2}y^3 + \eta_{3}y^2 + \eta_{4}y + \eta_{5}=0,
\end{eqnarray*}
where $\eta_{i}:= \eta_{i}(x, A, B, C, a, b, c)$, $i= 1, 2, 3, 4, 5$.

Computing the resultant of $f_{1}$ and $f_{2}$ with respect for $y$ and recalling that $A^{3}= A^{2}B-AB^{2}- C^{2}+C-1$, we derive
\begin{eqnarray*}
x^{4}(AB^2x^2 + (Cx - a)^2)^{4}(\alpha x +\beta)=0,
\end{eqnarray*}
where $\alpha := \alpha(A, B, C, a, b, c)$, $\beta := \beta(A, B, C, a, b, c)$. Since $a\neq0$, which means $x=0$ is not a solution of the above equation.
If $AB^2x^2 + (Cx - a)^2=0$, since $-A$ is a square, then replacing $A$ with $-A^2$ we have $x=\frac{a}{AB+C}$ or $x=-\frac{a}{AB-C}$.
Substituting $x=\frac{a}{AB+C}$ into Eqs.  (\ref{eq:q2.11}), (\ref{eq:q2.11'}), (\ref{eq:q2.11''}) eliminates $b$, $c$, we have
$$ A^{4}(A - 1)^{6}r_{1}(A,B)\cdot a^{16}=0, $$
where $r_{1}(A,B)=(A + B + 1)(A^2 + A + 1)^{6}(A^2 - AB - A + B^2 - B + 1)(A^3 + AB + 1)^{5}(A^4 - 2A^2B + AB^3 - A + B^2)(A^4 - A^2B + A - B^2)^{2}(A^4 + 2A^2B - A + B^2)$. Because $A\neq0, 1$ and $r_{1}(A,B)\neq0$, hence we have $a=0$, which is a contradiction. Similarly, we can discuss $x=-\frac{a}{AB-C}$ is also not a solution of the above equation.
Therefore Eq. (\ref{eq:q2.1}) has at most a solution $x=-\frac{\beta}{\alpha}$ and we complete the proof.
\qed

\subsection{A class of permutation pentanomials }
In the subsection, we give a class of permutation pentanomials of the form $f(x)= x^{q^2+q-1} + Ax^{q^2-q+1} + Bx^{q^2} + Cx^{q} + Dx\in \mathbb{F}_{q}[x]$.

\begin{theorem}\label{th:pp1}
Let $A, B, C, D\in\mathbb{F}_{q}$ with $C\neq0$, $B^2-4A=0$, $A^{3}-A^2C - ABD + AB + AC^2 - BC + D^2 - D + 1=0$, $AC^2-BCD+D^2\neq0$,
$A^3 - 2A^2C + AB + AC^2 - BC + 1\neq0$, $A^3 - AB^2C + AB + B^2C^2 - 2BC + 1\neq0$, $B^3+8\neq0$ and $AB - 3BC + 2C^3 + 2\neq0$, $m_{A, B, C, D}(t)=\lambda_{1}t^{3} + \lambda_{2}t^2 + \lambda_{3}t + \lambda_{4}\in \mathbb{F}_{q}[t]$ has no roots in $\mu_{q^2+q+1}$, $\lambda_{i}=\lambda_{i}(A, B, C, D), i=1, 2, 3, 4$, (see Appendix \ref{App}). Then
\begin{displaymath}\label{eq:pp1}
f(x)= x^{q^2+q-1} + Ax^{q^2-q+1} + Bx^{q^2} + Cx^{q} + Dx
\end{displaymath}
is a permutation pentanomial over $\mathbb{F}_{q^3}$.
\end{theorem}

\pf We prove that for each $a \in\mathbb{F}_{q^3}$, the equation
\begin{eqnarray}\label{eq:p1.1}
x^{q^2+q-1} + Ax^{q^2-q+1} + Bx^{q^2} + Cx^{q} + Dx = a
\end{eqnarray}
has at most a solution in $\mathbb{F}_{q^3}$.  We demonstrate $x=0$ if and only if $a=0$. Assume that $x\neq0$ is one solution of equation
\begin{eqnarray}\label{eq:p1.2}
x^{q^2+q-1} + Ax^{q^2-q+1} + Bx^{q^2} + Cx^{q} + Dx =0.
\end{eqnarray}
Setting $u=x^{q-1}$ implies $u^{q^2+q+1}=1$. It is necessary to prove that
\begin{eqnarray}\label{eq:p1.3}
u^{q+2} + Au^{q} + Bu^{q+1} + Cu + D =0
\end{eqnarray}
has no solution.
From Eq. (\ref{eq:p1.3}) we have $u^{q} = -\frac{Cu + D}{u^2+Bu+A}$ and $$u^{q^2} = \frac{C(Cu + D)(u^2+Bu+A)-D(u^2+Bu+A)^{2}}{(Cu+D)^{2} - B(Cu + D)(u^2+Bu+A)+A(u^2+Bu+A)^{2}}.$$
If $u^2+Bu+A=0$, then substituting it into Eq. (\ref{eq:p1.3}) gets $AC^{2}-BCD+D^2=0$, which is a contradiction.
If $(Cu+D)^{2} - B(Cu + D)(u^2+Bu+A)+A(u^2+Bu+A)^{2}=0$, then we obtain
\begin{eqnarray}\label{eq:p1d}
&& Au^4 + (2AB - BC)u^3 + (2A^2 + AB^2 - B^2C - BD + C^2)u^2 \nonumber \\
&+& (2A^2B - ABC - B^2D + 2CD)u + A^3 - ABD + D^2=0.
\end{eqnarray}
Raising Eq. (\ref{eq:p1d}) to the $q$-th power gets
\begin{eqnarray}\label{eq:p1dq}
&& (A^2C - AB - AC^2 + BC + D - 1)u^8 \nonumber \\
&+& (2A^2BC - 4AB^2 - 3ABC^2 + B^2CD + 4B^2C + 4BD - 4B - 2C^2D)u^7  \nonumber \\
&-&(2A^2BD + 4A^2B - 2A^2C^2 + 6AB^3 + 2AB^2C^2 - 5ABCD + 4AC^3 - 4AD + 4A - 3B^3CD   \nonumber \\
&-& 6B^3C + B^2C^3 - B^2D^2 - 6B^2D + 6B^2 + 7BC^2D - 4BC^2 - C^4 + 6CD^2 - 4CD + 4C)u^6  \nonumber \\
&-&(2A^2B^3C + 6A^2B^2D + 12A^2B^2 - A^2BC^2 - 4A^2CD + 4AB^4 - AB^3C^2 - 14AB^2CD   \nonumber \\
&-& 6AB^2C + 8ABC^3 - 12ABD + 12AB + 6AC^2D - 3B^4CD - 4B^4C + 2B^3C^3 - 3B^3D^2   \nonumber \\
&-& 4B^3D + 4B^3 + 10B^2C^2D - 6B^2C^2 - 3BC^4 + 14BCD^2 - 6BCD + 6BC - 2C^3D)u^5  \nonumber \\
&-& ( A^2B^4C + 6A^2B^3D + 12A^2B^3 + 2A^2B^2C^2 - 11A^2BCD + 6A^2BC + 2A^2C^3 - 2A^2D^2 \nonumber \\
&-& 6A^2D + 6A^2 + AB^5 - AB^4C^2 - 13AB^3CD - 12AB^3C + 4AB^2C^3 + 2AB^2D^2 \nonumber \\
&-& 12AB^2D + 6AB^2 + 10ABC^2D - 8ABC^2 + AC^4 + 12ACD^2 - 6ACD + 6AC \nonumber \\
&-& B^5CD - B^5C + B^4C^3 - 3B^4D^2 - B^4D + B^4 + 7B^3C^2D - 2B^2C^4 + 11B^2CD^2 \nonumber \\
&+& 6B^2CD + 6B^2C - 7BC^3D - 4BC^3 - 5BD^3 - 6B + 3C^2D^2 - 4C^2D + 4C^2)u^4  \nonumber \\
&-&(2A^2B^4D + 4A^2B^4 + A^2B^3C^2 - 11A^2B^2CD + 6A^2B^2C + 7A^2BC^3 - 4A^2BD^2 + 12A^2B  \nonumber \\
&-& 12A^2BD - 2A^2C^2D - 4AB^4CD - 6AB^4C + 4AB^3D^2 - 4AB^3D - 8AB^3 + 5AB^2C^2D  \nonumber \\
&-& 17AB^2C^2 - 2ABC^4 + 20ABCD^2 + 2ABCD + 6ABC - B^5D^2 + 2B^4C^2D + 2B^4C^2 \nonumber \\
&+& 4B^3CD^2 + 14B^3CD + 10B^3C - 5B^2C^3D - B^2C^3 - 10B^2D^3 - 12B^2 - 4BC^2D^2 \nonumber \\
&-& 9BC^2D + BC^2 + 8CD^3 - 8CD^2 + 8CD)u^3  \nonumber \\
&+& (4A^2B^3CD - 3A^2B^2C^3 + A^2B^2D^2 + 8A^2B^2D - 2A^2B^2 + 4A^2BC^2D + 2A^2BC^2 \nonumber \\
&-& 3A^2C^4 - 6A^2CD^2 + 8A^2CD - 8A^2C - 2AB^4D^2 + 6AB^4 + 8AB^3C^2 + 2AB^2C^4 \nonumber \\
&-& 9AB^2CD^2- 19AB^2CD - 4ABC^3D + 8ABC^3 + 6ABD^3 - 2ABD^2 - 6ABD \nonumber \\
&+& 8AB + 2AC^5 + 2AC^2D^2 - 2AC^2D + 2AC^2 - B^4CD^2 - 6B^4CD - 6B^4C \nonumber \\
&-& 2B^3C^3 + 5B^3D^3 + 6B^3 + 6B^2C^2D^2 +7B^2C^2D - 2B^2C^2 - 2BC^4 - 8BCD^3 \nonumber \\
&+& 17BCD^2 - 9BCD + 2C^3D^2 - 2C^3D + 2C^3 - 4D^4 + 4D^2 - 8D + 4)u^2  \nonumber \\
&-&(A^2B^3D^2 - 2A^2B^3D - 4A^2B^3 - 5A^2B^2C^2D - 5A^2B^2C^2 + 2A^2BC^4 + 8A^2BCD^2 \nonumber \\
&-& 2A^2BCD + 6A^2BC + 10AB^3CD + 2AB^3C + 5AB^2C^3D - 4AB^2C^3 - 6AB^2D^3 \nonumber \\
&+& 2AB^2D^2 + 6AB^2D - 8AB^2 - 3ABC^5 - 4ABC^2D^2 + 3ABC^2D - 5ABC^2 + 2AC^4D \nonumber \\
&-& 4ACD^2 + 4ACD - 2B^3C^2D + 2B^3C^2 + 3B^2C^4 + B^2CD^3 - 8B^2CD^2 + 2B^2C \nonumber \\
&-& 3BC^3D^2 + BC^3D - 3BC^3 + 4BD^4 - 4BD^2 + 8BD - 4B + 2C^2D^3 - 2C^2D^2 + 2C^2D)u  \nonumber \\
&-& A^2B^2CD^2 - 2A^2B^2CD + A^2B^2C + A^2BC^3D + 3A^2BC^3 + A^2BD^3 - 2A^2BD^2 - A^2BD \nonumber \\
&+& 2A^2B - 2A^2C^2D^2 + A^2C^2D - A^2C^2 + AB^3D^2 - AB^3 + 2AB^2C^2D^2 - 2AB^2C^2D \nonumber \\
&-& 3AB^2C^2 - 3ABC^4D - ABC^4 + 6ABCD^2 - 4ABCD + 2ABC + AC^6 + 2AC^3D^2 \nonumber \\
&-& 3AC^3D + 3AC^3 - AD^4 + AD^3 - 2AD + A + B^3CD + B^3C + 2B^2C^3D + B^2C^3 \nonumber \\
&-& B^2D^3 - B^2 - BC^5 - 2BC^2D^3 + BC^2D - 3BC^2 + C^4D^2 - C^4D + C^4 + CD^4 \nonumber \\
&-& 3CD^3 + 5CD^2 - 4CD + 2C=0.
\end{eqnarray}
By Eqs. (\ref{eq:p1d}) and (\ref{eq:p1dq}) eliminates $u$ we have
\begin{eqnarray*}
(A^3 - 2A^2C + AB + AC^2 - BC + 1)^{12}(A^3 - AB^2C + AB + B^2C^2 - 2BC + 1)^{12}=0,
\end{eqnarray*}
since $A^3 - 2A^2C + AB + AC^2 - BC + 1\neq0$ and $A^3 - AB^2C + AB + B^2C^2 - 2BC + 1\neq0$, which means  Eq. (\ref{eq:p1d}) has no solution, thereby $(Cu+D)^{2} - B(Cu + D)(u^2+Bu+A)+A(u^2+Bu+A)^{2}\neq0$.

Note that $u^{q^2+q+1}=1$, we obtain
\begin{eqnarray*}
&& (u^2+Bu+A) \big[ (A - CD)u^4 + (2AB - BCD - BC + C^3 - D^2)u^3 \nonumber \\
&+& (2A^2 + AB^2 - ACD - B^2C - BD^2 - BD + 2C^2D + C^2)u^2 \nonumber \\
&+& (2A^2B - ABC - AD^2 - B^2D + CD^2 + 2CD)u + A^3 - ABD + D^2 \big]=0.
\end{eqnarray*}
Since $AC^{2}-BCD+D^2\neq0$, we get $u^2+Bu+A=0$ has no solution and using $A^{3}= A^2C + ABD - AB - AC^2 + BC - D^2 + D - 1$ we deduce
\begin{eqnarray*}
&& (A - CD)u^4 + (2AB - BCD - BC + C^3 - D^2)u^3 \nonumber \\
&+& (2A^2 + AB^2 - ACD - B^2C - BD^2 - BD + 2C^2D + C^2)u^2 \nonumber \\
&+& (2A^2B - ABC - AD^2 - B^2D + CD^2 + 2CD)u + A^2C - AB - AC^2 + BC + D - 1 =0.
\end{eqnarray*}
Raising both sides of the above equation to the $q$-th power leads to
\begin{eqnarray}\label{eq:pp2.5}
&&  (A^2C - AB - AC^2 + BC + D - 1)u^8 \nonumber \\
&+& (2A^2BC - 4AB^2 - 3ABC^2 + ACD^2 + B^2CD + 4B^2C + 4BD - 4B - C^2D^2 - 2C^2D)u^7   \nonumber \\
&-& (2A^2BD + 4A^2B - 2A^2C^2 + 6AB^3 + 2AB^2C^2 - 3ABCD^2 - 5ABCD + AC^3D + 4AC^3 \nonumber \\
&-& AD^3 - 4AD + 4A - 3B^3CD - 6B^3C + B^2C^3 - B^2D^2 - 6B^2D + 6B^2 + 4BC^2D^2 \nonumber \\
&+& 7BC^2D - 4BC^2 - 2C^4D - C^4 + CD^3 + 6CD^2 - 4CD + 4C)u^6 \nonumber \\
&-& (2A^2B^3C + 6A^2B^2D + 12A^2B^2 - A^2BC^2 - 3A^2CD^2 - 4A^2CD + 4AB^4 - AB^3C^2 \nonumber \\
&-& 3AB^2CD^2 - 14AB^2CD - 6AB^2C + 2ABC^3D + 8ABC^3 - 3ABD^3 - 12ABD + 12AB \nonumber \\
&+& 5AC^2D^2 + 6AC^2D - 3B^4CD - 4B^4C + 2B^3C^3 - 3B^3D^2 - 4B^3D + 4B^3 + 5B^2C^2D^2 \nonumber \\
&+& 10B^2C^2D - 6B^2C^2 - 5BC^4D - 3BC^4 + 5BCD^3 + 14BCD^2 - 6BCD + 6BC + C^6 \nonumber \\
&-& 5C^3D^2 - 2C^3D)u^5 -(A^2B^4C + 6A^2B^3D + 12A^2B^3 + 2A^2B^2C^2 - 6A^2BCD^2 \nonumber \\
&-& 11A^2BCD + 6A^2BC + 2A^2C^3D + 2A^2C^3 - 3A^2D^3 - 2A^2D^2 - 6A^2D + 6A^2 + AB^5 \nonumber \\
&-& AB^4C^2 - AB^3CD^2 - 13AB^3CD - 12AB^3C + AB^2C^3D + 4AB^2C^3 - 3AB^2D^3 \nonumber \\
&+& 2AB^2D^2 - 12AB^2D + 6AB^2 + 12ABC^2D^2 + 10ABC^2D - 8ABC^2 - 4AC^4D + AC^4 \nonumber \\
&+& 4ACD^3 + 12ACD^2 - 6ACD + 6AC - B^5CD - B^5C + B^4C^3 - 3B^4D^2 - B^4D + B^4 \nonumber \\
&+& 2B^3C^2D^2 + 7B^3C^2D - 3B^2C^4D - 2B^2C^4 + 7B^2CD^3 + 11B^2CD^2 + 6B^2CD + 6B^2C \nonumber \\
&+& BC^6 - 12BC^3D^2 - 7BC^3D - 4BC^3 + BD^4 - 5BD^3 - 6B + 4C^5D - 5C^2D^3 + 3C^2D^2 \nonumber \\
&-& 4C^2D + 4C^2)u^4-(2A^2B^4D + 4A^2B^4 + A^2B^3C^2 - 3A^2B^2CD^2 - 11A^2B^2CD + 6A^2B^2C \nonumber \\
&+& 2A^2BC^3D + 7A^2BC^3 - 6A^2BD^3 - 4A^2BD^2 - 12A^2BD + 12A^2B + 4A^2C^2D^2 \nonumber \\
&-& 2A^2C^2D - 4AB^4CD - 6AB^4C - AB^3D^3 + 4AB^3D^2 - 4AB^3D - 8AB^3 + 7AB^2C^2D^2 \nonumber \\
&+& 5AB^2C^2D - 17AB^2C^2 - 5ABC^4D - 2ABC^4 + 9ABCD^3 + 23ABCD^2 + 2ABCD \nonumber \\
&+& 6ABC + AC^6 - 6AC^3D^2 - B^5D^2 + 2B^4C^2D + 2B^4C^2 + 3B^3CD^3 + 4B^3CD^2 \nonumber \\
&+& 14B^3CD + 10B^3C - 7B^2C^3D^2 - 5B^2C^3D - B^2C^3 + 2B^2D^4 - 10B^2D^3 - 12B^2 \nonumber \\
&+& 3BC^5D - 10BC^2D^3 - 7BC^2D^2 - 9BC^2D + BC^2 + 7C^4D^2 + 5CD^3 - 5CD^2 + 8CD)u^3 \nonumber \\
&+& (4A^2B^3CD - 3A^2B^2C^3 + 3A^2B^2D^3 + A^2B^2D^2 + 8A^2B^2D - 2A^2B^2 - 5A^2BC^2D^2 \nonumber \\
&+& 4A^2BC^2D + 2A^2BC^2 + A^2C^4D - 3A^2C^4 - 2A^2CD^3 - 6A^2CD^2 + 8A^2CD - 8A^2C \nonumber \\
&-& 2AB^4D^2 + 6AB^4 + 8AB^3C^2 + 2AB^2C^4 - 5AB^2CD^3 - 12AB^2CD^2 - 19AB^2CD \nonumber \\
&+& 7ABC^3D^2 - 3ABC^3D + 8ABC^3 + ABD^4 + 3ABD^3 - 2ABD^2 - 6ABD + 8AB \nonumber \\
&-& 2AC^5D + 2AC^5 + 4AC^2D^3 + 2AC^2D^2 - 2AC^2D + 2AC^2 - B^4CD^2 - 6B^4CD \nonumber \\
&-& 6B^4C - 2B^3C^3 - B^3D^4 + 5B^3D^3 + 6B^3 + 5B^2C^2D^3 + 9B^2C^2D^2 + 7B^2C^2D \nonumber \\
&-& 2B^2C^2 - 3BC^4D^2 - BC^4D - 2BC^4 + BCD^4 - 2BCD^3 + 14BCD^2 - 9BCD \nonumber \\
&-& 6C^3D^3 + C^3D^2 - C^3D + 2C^3 - 2D^5 - D^4 - 3D^3 + 4D^2 - 8D + 4)u^2 \nonumber \\
&-&(A^2B^3D^2 - 2A^2B^3D - 4A^2B^3 - 5A^2B^2C^2D - 5A^2B^2C^2 + 2A^2BC^4 + 3A^2BCD^3 \nonumber \\
&+& 9A^2BCD^2 - 2A^2BCD + 6A^2BC - A^2C^3D^2 + 10AB^3CD + 2AB^3C + 5AB^2C^3D \nonumber \\
&-& 4AB^2C^3 - AB^2D^4 - 3AB^2D^3 + 2AB^2D^2 + 6AB^2D - 8AB^2 - 3ABC^5 - 2ABC^2D^3 \nonumber \\
&-& 7ABC^2D^2 + 3ABC^2D - 5ABC^2 + AC^4D^2 + 2AC^4D - 2ACD^4 - ACD^3 - 3ACD^2 \nonumber \\
&+& 4ACD - 2B^3C^2D + 2B^3C^2 + 3B^2C^4 - B^2CD^4 - 2B^2CD^3 - 8B^2CD^2 + 2B^2C \nonumber \\
&+& BC^3D^3 - BC^3D^2 + BC^3D - 3BC^3 + 2BD^5 + BD^4 + 3BD^3 - 4BD^2 + 8BD - 4B \nonumber \\
&+& 2C^2D^4 + 4C^2D^3 - 4C^2D^2 + 2C^2D)u \nonumber \\
&-& A^2B^2CD^2 - 2A^2B^2CD + A^2B^2C + A^2BC^3D + 3A^2BC^3 - 2A^2BD^2 - A^2BD + 2A^2B \nonumber \\
&-& 2A^2C^2D^2 + A^2C^2D - A^2C^2 + AB^3D^2 - AB^3 + 2AB^2C^2D^2 - 2AB^2C^2D - 3AB^2C^2 \nonumber \\
&-& 3ABC^4D - ABC^4 + 2ABCD^3 + 6ABCD^2 - 4ABCD + 2ABC + AC^6 + 2AC^3D^2 \nonumber \\
&-& 3AC^3D + 3AC^3 - 2AD + A + B^3CD + B^3C + 2B^2C^3D + B^2C^3 - B^2D^3 - B^2 - BC^5 \nonumber \\
&-& 3BC^2D^3 + BC^2D - 3BC^2 + C^4D^2 - C^4D + C^4 - 2CD^3 + 5CD^2 - 4CD + 2C  = 0.
 \end{eqnarray}
Through a series of computations we obtain $m_{A, B, C, D}(u)=\lambda_{1}u^{3} + \lambda_{2}u^2 + \lambda_{3}u + \lambda_{4}=0 $, $\lambda_{i}=\lambda_{i}(A, B, C, D), i=1, 2, 3, 4$, (see Appendix \ref{App}), which contradicts $m_{A, B, C, D}(u)=0 $ has no solution in $\mu_{q^2+q+1}$ and therefore $f(x)=0$ has a unique solution $x=0$.

If $a\neq0$, we demonstrate that Eq. (\ref{eq:p1.1}) has one nonzero solution. Let $y=x^{q}$, $z=y^{q}$, $b=a^{q}$ and $c=b^{q}$, then we obtain the system of equations
\begin{numcases}{}
y^2z + Ax^2z + Bxyz + Cxy^2 + Dx^2y -axy =0,  \label{eq:p1.11}\\
z^2x + Ay^2x + Bxyz + Cyz^2 + Dy^2z -byz =0,  \label{eq:p1.11'}\\
x^2y + Az^2y + Bxyz + Czx^2 + Dz^2x -czx  =0. \label{eq:p1.11''}
\end{numcases}
Eliminating the indeterminate $z$ by (\ref{eq:p1.11}) and (\ref{eq:p1.11'}), we have
\begin{eqnarray*}
f_{1} \triangleq \xi_{1}y^4 + \xi_{2}y^3 + \xi_{3}y^2 + \xi_{4}y + \xi_{5}=0,
\end{eqnarray*}
where $\xi_{i}:= \xi_{i}(x, A, B, C, D, a, b, c)$, $i= 1, 2, 3, 4, 5$.

Furthermore, by Eqs. (\ref{eq:p1.11}) and (\ref{eq:p1.11''}) we get
\begin{eqnarray*}
f_{2} \triangleq \eta_{1}y^4 + \eta_{2}y^3 + \eta_{3}y^2 + \eta_{4}y + \eta_{5}=0,
\end{eqnarray*}
where $\eta_{i}:= \eta_{i}(x, A, B, C, D, a, b, c)$, $i= 1, 2, 3, 4, 5$.

Computing the resultant of $f_{1}$ and $f_{2}$ with respect for $y$ and recalling that $A^{3}= A^2C + ABD - AB - AC^2 + BC - D^2 + D - 1$, we have
\begin{eqnarray*}
x^{4}((AC^2 - BCD + D^2)x^2 + (BC - 2D)ax + a^2)^{4}(\alpha x +\beta)=0,
\end{eqnarray*}
where $\alpha := \alpha(A, B, C, D, a, b, c)$, $\beta := \beta(A, B, C, D, a, b, c)$. Since $a\neq0$, which means $x=0$ is not a solution of the above equation.
If $(AC^2 - BCD + D^2)x^2 + (BC - 2D)ax + a^2=0$, since $B^2-4A=0$,
then we obtain $x=\frac{2a(2D-BC)}{B^2C^2-4BCD+4D^2}$, which implies
$y=\frac{2b(2D-BC)}{B^2C^2-4BCD+4D^2}$ and $z=\frac{2c(2D-BC)}{B^2C^2-4BCD+4D^2}$.
Substituting them into Eqs.  (\ref{eq:p1.11}) and factoring the equation, we have
$$(BC - 2D)(Ba + 2b)(Bac + 2Cab + 2bc)=0,$$
if $b=-\frac{B}{2}a$, then we deduce $(B^3+8)a=0$, however $B^3+8\neq0$, a contradiction. Since $BC - 2D \neq0$, we obtain
$$Bac + 2Cab + 2bc=0,$$
raising the equation to the $q$-th and $q^2$-th power and eliminating $b$, $c$ with $B^2=4A$, we deduce
$$ 8C(AB - 3BC + 2C^3 + 2)a=0.$$
Because $C\neq0$ and $AB - 3BC + 2C^3 + 2\neq0$, hence we have $a=0$, which is a contradiction.
Therefore Eq. (\ref{eq:p1.1}) has at most a solution $x=-\frac{\beta}{\alpha}$ and we complete the proof.
\qed

Despite that there are many conditions on these coefficients in several of the previous theorems, these conditions can be easily checked by using computer programs. In Table \ref{Table2}, we provide some explicit classes of permutation polynomials as examples to demonstrate the previous theorems.

\begin{table}[!htbp]
\caption{Permutation polynomials over $\mathbb{F}_{q^3}$ }\label{Table2}
\begin{tabular}{ccc}
\toprule
      PPs   &                                     Conditions &                              ~~~~~~~  Theorems  \\
	\midrule
$x^{q^2+q-1} +3x$  &                         $q=7^k$, $k$ is a integer  &                   ~~~~~~~~~~  Theorem \ref{th:pb}  \\

$x^{q^2-q+1} + 2x^{q^3-q^2+q} +3x$  &        $q=5^k$, $k$ is a integer  &                 ~~~~~~~~~~   Theorem \ref{th:pt1}  \\

$x^{q^2+q-1} + 2x^{q^2} +4x$  &              $q=13^k$, $k\not\equiv 2~({\rm mod}\ 3)$  &   ~~~~~~~~~~  Theorem \ref{th:pt2}        \\
$x^{q^2+q-1} + 2x^{q} + 3x$  &                $q=7^k$, $k\not\equiv 2~({\rm mod}\ 3)$  &   ~~~~~~~~~~  Theorem \ref{th:pt3}    \\
$x^{q^2+q-1} + x^{q^2}+ 2x^{q} + 3x$  &      $q=5^k$, $k\not\equiv 2~({\rm mod}\ 3)$  &   ~~~~~~~~~~  Theorem \ref{th:pq1}    \\
$x^{q^2+q-1} + x^{q^2-q+1}+ 4x^{q^2} + 2x$  & $q=5^k$, $k\not\equiv 2~({\rm mod}\ 3)$   &   ~~~~~~~~~~  Theorem \ref{th:pq3}    \\
$x^{q^2+q-1} + x^{q^3-q^2+q}+ 2x^{q} + 6x$  & $q=11^k$, $k\not\equiv 2~({\rm mod}\ 3)$  &  ~~~~~~~~~~  Theorem \ref{th:pq4}    \\
$x^{q^2+q-1} + 3x^{q^2-q+1}+ x^{q} + x$  &   $q=11^k$, $k$ is a integer              &   ~~~~~~~~~~ Theorem \ref{th:pq2}    \\

$x^{q^2+q-1} + 4x^{q^2-q+1}+ x^{q^2}+ 2x^{q} + 3x$  &   $q=5^k$, $k\not\equiv 2~({\rm mod}\ 3)$  &   ~~~~~~~~~~  Theorem \ref{th:pp1}    \\

\bottomrule
\end{tabular}
\end{table}

\section{Conclusion}\label{conclu}
By using the multivariate method, we construct a class of complete permutation binomials with the coefficient over $\mathbb{F}_{q^3}$, several classes of permutation trinomials, permutation quadrinomials and one class of permutation pentanomials in terms of their coefficients in $\mathbb{F}_{q}$. Their permutation properties are proved by using the resultant elimination method which is a useful tool in this paper.

\section{Appendix}\label{App}

In this section, we list some necessary MAGMA programs used in the proofs of the results in the previous sections and some equations.

Firstly, we introduce a ``Substitution" function {\rm\cite{Bartoli2018}} which will be applied over and over again.

{\scriptsize
\begin{verbatim}
Substitution := function (pol, m, p)
	e := 0;
	New := R! pol;
	while e eq 0 do
		N := R!0;
		T := Terms(New);
		i:= 0;
		for t in T do
			if IsDivisibleBy(t,m) eq true then
				Q := R! (t/m);
				i := 1;
				N := R!(N + Q* p);
			else
				N := R!(N + t);
			end if;
		end for;
		if i eq 0 then
			return New;
		else	
			New := R!N;
		end if;	
	end while;
end function;
\end{verbatim}
}

\subsection{Theorem \ref{th:pb}}

{\scriptsize
\begin{verbatim}
// The case of a ne 0:
R<x,y,z,A,B,C,a,b,c> := PolynomialRing(Integers(),9);
f1 := y*z + A*x^2 - a*x;
f2 := z*x + B*y^2 - b*y;
f3 := x*y + C*z^2 - c*z;
R1 := R!(Resultant(f1,f2,z));
R2 := R!(Resultant(f2,f3,z));
RR := R!(Resultant(R1,R2,y));
Factorization(RR);
\end{verbatim}
}

\subsection{Theorem  \ref{th:pt1}}

{\scriptsize
\begin{verbatim}
// The case of a ne 0:
R<x,y,z,A,B,a,b,c> := PolynomialRing(Integers(),8);
f1 := x*z^2 + A*x*y^2 + B*x*y*z - a*z*y;
f2 := Evaluate(f1,[y,z,x,A,B,b,c,a]);
f3 := Evaluate(f2,[y,z,x,A,B,b,c,a]);
R1 := R!(Resultant(f1,f3,z)/x/y^2);
R2 := R!(Resultant(f2,f3,z)/x^2/y);
RR := R!(Resultant(R1,R2,y)/c^8/x^4);
RR := Substitution(RR, A^3, -A*B-1);
Factorization(RR);
\end{verbatim}
}

\subsection{ Theorem \ref{th:pt2}}

{\scriptsize
\begin{verbatim}
// The case of a ne 0:
R<x,y,z,A,C,a,b,c> := PolynomialRing(Integers(),8);
f1 := y*z + A*z*x + C*x^2 - a*x;
f2 := Evaluate(f1,[y,z,x,A,C,b,c,a]);
f3 := Evaluate(f2,[y,z,x,A,C,b,c,a]);
R1 := R!(Resultant(f1,f2,z));
R2 := R!(Resultant(f1,f3,z)/x);
RR := R!(Resultant(R1,R2,y)/x^2/(C*x-a)^3);
RR := Substitution(RR, C^2, C-1);
Factorization(RR);
\end{verbatim}
}

\subsection{ Theorem \ref{th:pt3} }

{\scriptsize
\begin{verbatim}
// The case of a ne 0:
R<x,y,z,B,C,a,b,c> := PolynomialRing(Integers(),8);
f1 := y*z + B*y*x + C*x^2 - a*x;
f2 := Evaluate(f1,[y,z,x,B,C,b,c,a]);
f3 := Evaluate(f2,[y,z,x,B,C,b,c,a]);
R1 := R!(Resultant(f1,f2,z));
R2 := R!(Resultant(f1,f3,z)/x);
RR := R!(Resultant(R1,R2,y)/x^2/(C*x-a)^3);
RR := Substitution(RR, C^2, C-1);
Factorization(RR);
\end{verbatim}
}

\subsection{ Theorem \ref{th:pq1}}

{\scriptsize
\begin{verbatim}
// The case of a ne 0:
R<x,y,z,A,B,C,a,b,c> := PolynomialRing(Integers(),9);
f1 := y*z + A*x*z + B*x*y + C*x^2 - a*x;
f2 := Evaluate(f1,[y,z,x,A,B,C,b,c,a]);
f3 := Evaluate(f2,[y,z,x,A,B,C,b,c,a]);
R1 := R!(Resultant(f1,f2,z));
R2 := R!(Resultant(f1,f3,z)/x);
RR := R!(Resultant(R1,R2,y));
RR := Substitution(RR, A*B, C^2-C+1);
Factorization(RR);
\end{verbatim}
}

\subsection{ Theorem \ref{th:pq3}}

{\scriptsize
\begin{verbatim}
// The case of a ne 0:
R<x,y,z,A,B,C,a,b,c> := PolynomialRing(Integers(),9);
f1 := y^2*z + A*x^2*z + B*x*y*z + C*x^2*y - a*x*y;
f2 := Evaluate(f1,[y,z,x,A,B,C,b,c,a]);
f3 := Evaluate(f2,[y,z,x,A,B,C,b,c,a]);
R1 := R!(Resultant(f1,f2,z)/x/y^2);
R2 := R!(Resultant(f1,f3,z)/x^2/y);
RR := R!(Resultant(R1,R2,y));
RR := Substitution(RR, A^3, A*B*C-A*B-C^2+C-1);
Factorization(RR);
\end{verbatim}
}

\subsection{ Theorem \ref{th:pq4}}

{\scriptsize
\begin{verbatim}
// The case of a ne 0:
R<x,y,z,A,B,C,a,b,c> := PolynomialRing(Integers(),9);
f1 := y*z^2 + A*x^2*y + B*x*y*z + C*x^2*z - a*x*z;
f2 := Evaluate(f1,[y,z,x,A,B,C,b,c,a]);
f3 := Evaluate(f2,[y,z,x,A,B,C,b,c,a]);
R1 := R!(Resultant(f1,f2,z)/x^2/y);
R2 := R!(Resultant(f2,f3,z)/x/y^2);
RR := R!(Resultant(R1,R2,y)/x^4/(x^2*C^2 + x*B*C*b + A*b^2 )^4);
RR := Substitution(RR, A^3, A*B*C-A*B-C^2+C-1);
Factorization(RR);
\end{verbatim}
}

\subsection{ Theorem \ref{th:pq2}}

{\scriptsize
\begin{verbatim}
// The case of a ne 0:
R<x,y,z,A,B,C,a,b,c> := PolynomialRing(Integers(),9);
f1 := y^2*z + A*x^2*z + B*x*y^2 + C*x^2*y - a*x*y;
f2 := Evaluate(f1,[y,z,x,A,B,C,b,c,a]);
f3 := Evaluate(f2,[y,z,x,A,B,C,b,c,a]);
R1 := R!(Resultant(f1,f2,z)/x/y^2);
R2 := R!(Resultant(f1,f3,z)/x^2/y);
RR := R!(Resultant(R1,R2,y)/x^4/((A*B^2*x^2+(C*x-a)^2)^4));
RR := Substitution(RR, A^3, A^2*B-A*B^2-C^2+C-1);
Factorization(RR);
\end{verbatim}
}

\begin{scriptsize}
$m_{A,B,C}(u)= (A^2B^{16} + 6A^2B^{13}C^2 - 6A^2B^{13}C + 6A^2B^{13} + 9A^2B^{10}C^4 -
        30A^2B^{10}C^3 + 24A^2B^{10}C^2 - 30A^2B^{10}C + 9A^2B^{10} +
        4A^2B^7C^6 - 30A^2B^7C^5 + 30A^2B^7C^4 - 32A^2B^7C^3 +
        30A^2B^7C^2 - 30A^2B^7C + 4A^2B^7 - 15A^2B^4C^6 +
        70A^2B^4C^5 - 85A^2B^4C^4 + 70A^2B^4C^3 - 15A^2B^4C^2 -
        6A^2BC^5 - AB^{17} - 3AB^{14}C^2 + 9AB^{14}C - 3AB^{14} +
        30AB^{11}C^3 - 30AB^{11}C^2 + 30AB^{11}C + 2AB^8C^6 + 24AB^8C^5
        - 78AB^8C^4 + 86AB^8C^3 - 78AB^8C^2 + 24AB^8C + 2AB^8 +
        6AB^5C^7 - 48AB^5C^6 + 76AB^5C^5 - 97AB^5C^4 + 76AB^5C^3
        - 48AB^5C^2 + 6AB^5C + 15AB^2C^6 - 15AB^2C^5 + 15AB^2C^4
        - B^{15}C^2 + B^{15}C - B^{15} - 3B^{12}C^4 + 12B^{12}C^3 - 15B^{12}C^2 +
        12B^{12}C - 3B^{12} - 2B^9C^6 + 30B^9C^5 - 60B^9C^4 + 86B^9C^3 -
        60B^9C^2 + 30B^9C - 2B^9 - B^6C^8 + 22B^6C^7 - 79B^6C^6 +
        134B^6C^5 - 170B^6C^4 + 134B^6C^3 - 79B^6C^2 + 22B^6C - B^6 -
        20B^3C^7 + 55B^3C^6 - 75B^3C^5 + 55B^3C^4 - 20B^3C^3 - C^6 )u^3
+ (2A^2B^{15}C - A^2B^{15} + 9A^2B^{12}C^3 - 15A^2B^{12}C^2 + 20A^2B^{12}C
        - A^2B^{12} + 11A^2B^9C^5 - 51A^2B^9C^4 + 81A^2B^9C^3 -
        89A^2B^9C^2 + 39A^2B^9C + A^2B^9 + 5A^2B^6C^7 - 46A^2B^6C^6
        + 103A^2B^6C^5 - 171A^2B^6C^4 + 167A^2B^6C^3 - 121A^2B^6C^2
        + 28A^2B^6C - A^2B^6 - 10A^2B^3C^7 + 50A^2B^3C^6 -
        75A^2B^3C^5 + 70A^2B^3C^4 - 20A^2B^3C^3 - A^2C^6 - AB^{16}C -
        AB^{16} + 2AB^{13}C^2 + 3AB^{13}C - 8AB^{13} + 7AB^{10}C^5 +
        9AB^{10}C^4 + 5AB^{10}C^2 + 31AB^{10}C - 12AB^{10} + 5AB^7C^7 +
        6AB^7C^6 - 48AB^7C^5 + 76AB^7C^4 - 62AB^7C^3 + 6AB^7C^2 +
        27AB^7C - 4AB^7 + 5AB^4C^8 - 45AB^4C^7 + 100AB^4C^6 -
        155AB^4C^5 + 135AB^4C^4 - 85AB^4C^3 + 15AB^4C^2 + 5ABC^7
        - 6ABC^6 + 6ABC^5 - B^{14}C^3 - B^{14} - B^{11}C^5 + 6B^{11}C^4 -
        6B^{11}C^3 - B^{11}C^2 + 6B^{11}C - 6B^{11} + B^8C^7 + 21B^8C^6 -
        51B^8C^5 + 87B^8C^4 - 57B^8C^3 + 21B^8C^2 + 14B^8C - 6B^8 -
        B^5C^9 + 24B^5C^8 - 105B^5C^7 + 231B^5C^6 - 339B^5C^5 +
        321B^5C^4 - 200B^5C^3 + 69B^5C^2 - 6B^5C - 10B^2C^8 +
        30B^2C^7 - 45B^2C^6 + 36B^2C^5 - 15B^2C^4 )u^2
+(C+1)(A^2B^{14}C + 3A^2B^{11}C^3 - 6A^2B^{11}C^2 + 6A^2B^{11}C + 2A^2B^{11} +
        2A^2B^8C^5 - 21A^2B^8C^4 + 30A^2B^8C^3 - 35A^2B^8C^2 +
        6A^2B^8C + 3A^2B^8 + A^2B^5C^7 - 18A^2B^5C^6 + 51A^2B^5C^5
        - 82A^2B^5C^4 + 81A^2B^5C^3 - 48A^2B^5C^2 + 6A^2B^5C +
        10A^2B^2C^6 - 15A^2B^2C^5 + 15A^2B^2C^4 - AB^{15} + 3AB^{12}C^3
        - 3AB^{12}C^2 + 8AB^{12}C - 5AB^{12} + 7AB^9C^5 - 9AB^9C^4 +
        21AB^9C^3 - 19AB^9C^2 + 27AB^9C - 5AB^9 + 3AB^6C^7 -
        10AB^6C^6 + AB^6C^5 - 7AB^6C^4 + 5AB^6C^3 - 25AB^6C^2 +
        16AB^6C - AB^6 - 10AB^3C^7 + 30AB^3C^6 - 45AB^3C^5 +
        40AB^3C^4 - 20AB^3C^3 - AC^6 - B^{13}C^2 + B^{13}C - B^{13} +
        2B^{10}C^5 - 3B^{10}C^4 + 9B^{10}C^3 - 11B^{10}C^2 + 10B^{10}C - 4B^{10}
        + 3B^7C^7 - 2B^7C^6 + 2B^7C^5 + B^7C^4 + 14B^7C^3 - 17B^7C^2
        + 17B^7C - 3B^7 + 5B^4C^8 - 35B^4C^7 + 80B^4C^6 - 125B^4C^5 +
        110B^4C^4 - 65B^4C^3 + 15B^4C^2 + 5BC^7 - 5BC^6 + 6BC^5)u
- A^2B^{16} - 6A^2B^{13}C^2 +
        6A^2B^{13}C - 6A^2B^{13} - 9A^2B^{10}C^4 + 30A^2B^{10}C^3 -
        24A^2B^{10}C^2 + 30A^2B^{10}C - 9A^2B^{10} - 4A^2B^7C^6 +
        30A^2B^7C^5 - 30A^2B^7C^4 + 32A^2B^7C^3 - 30A^2B^7C^2 +
        30A^2B^7C - 4A^2B^7 + 15A^2B^4C^6 - 70A^2B^4C^5 +
        85A^2B^4C^4 - 70A^2B^4C^3 + 15A^2B^4C^2 + 6A^2BC^5 + AB^17
        + 3AB^{14}C^2 - 9AB^{14}C + 3AB^{14} - 30AB^{11}C^3 + 30AB^{11}C^2 -
        30AB^{11}C - 2AB^8C^6 - 24AB^8C^5 + 78AB^8C^4 - 86AB^8C^3 +
        78AB^8C^2 - 24AB^8C - 2AB^8 - 6AB^5C^7 + 48AB^5C^6 -
        76AB^5C^5 + 97AB^5C^4 - 76AB^5C^3 + 48AB^5C^2 - 6AB^5C -
        15AB^2C^6 + 15AB^2C^5 - 15AB^2C^4 + B^{15}C^2 - B^{15}C + B^{15} +
        3B^{12}C^4 - 12B^{12}C^3 + 15B^{12}C^2 - 12B^{12}C + 3B^{12} + 2B^9C^6
        - 30B^9C^5 + 60B^9C^4 - 86B^9C^3 + 60B^9C^2 - 30B^9C + 2B^9 +
        B^6C^8 - 22B^6C^7 + 79B^6C^6 - 134B^6C^5 + 170B^6C^4 -
        134B^6C^3 + 79B^6C^2 - 22B^6C + B^6 + 20B^3C^7 - 55B^3C^6 +
        75B^3C^5 - 55B^3C^4 + 20B^3C^3 + C^6 $.
\end{scriptsize}

\subsection{ Theorem \ref{th:pp1}}

{\scriptsize
\begin{verbatim}
// The case of a ne 0:
R<x,y,z,A,B,C,D,a,b,c> := PolynomialRing(Integers(),10);
f1 := y^2*z + A*x^2*z + B*x*y*z + C*x*y^2 + D*x^2*y - a*x*y;
f2 := Evaluate(f1,[y,z,x,A,B,C,D,b,c,a]);
f3 := Evaluate(f2,[y,z,x,A,B,C,D,b,c,a]);
R1 := R!(Resultant(f1,f2,z)/x/y^2);
R2 := R!(Resultant(f1,f3,z)/x^2/y);
RR := R!(Resultant(R1,R2,y)/x^4/((A*C^2 - B*C*D + D^2)*x^2 + (B*C*a - 2*D*a)*x + a^2)^4);
RR := Substitution(RR, A^3, A^2*C + A*B*D - A*B - A*C^2 + B*C - D^2 + D - 1);
Factorization(RR);
\end{verbatim}
}

\begin{scriptsize}
$ m_{A,B,C,D}(u)= (A^2B^5C^6D^3 + 2A^2B^5C^6D^2 + 3A^2B^5C^6D - 2A^2B^5C^6 +
        3A^2B^4C^8D^2 - 8A^2B^4C^8D - 11A^2B^4C^8 - 6A^2B^4C^5D^3
        - 12A^2B^4C^5D^2 + 12A^2B^4C^5D - 4A^2B^4C^5 -
        6A^2B^3C^{10}D^2 - 27A^2B^3C^{10}D - 6A^2B^3C^{10} -
        3A^2B^3C^7D^4 - 10A^2B^3C^7D^3 + 61A^2B^3C^7D^2 +
        20A^2B^3C^7D + 5A^2B^3C^7 + 15A^2B^3C^4D^3 -
        30A^2B^3C^4D^2 + 15A^2B^3C^4D - 2A^2B^2C^{12}D^2 +
        15A^2B^2C^{12} + 22A^2B^2C^9D^3 + 66A^2B^2C^9D^2 -
        24A^2B^2C^9D + 38A^2B^2C^9 + A^2B^2C^6D^6 + 2A^2B^2C^6D^5
        + 21A^2B^2C^6D^4 - 130A^2B^2C^6D^3 + 28A^2B^2C^6D^2 -
        60A^2B^2C^6D + 15A^2B^2C^6 + 40A^2B^2C^3D^3 -
        20A^2B^2C^3D^2 + 3A^2BC^{14}D + 4A^2BC^{14} + 8A^2BC^{11}D^3 -
        12A^2BC^{11}D^2 - 30A^2BC^{11}D + 3A^2BC^{11} + 3A^2BC^8D^5 -
        36A^2BC^8D^4 - 47A^2BC^8D^3 + 66A^2BC^8D^2 - 69A^2BC^8D
        - 4A^2BC^8 - 6A^2BC^5D^6 - 12A^2BC^5D^5 + 131A^2BC^5D^4 -
        122A^2BC^5D^3 + 135A^2BC^5D^2 - 36A^2BC^5D + A^2BC^5 -
        30A^2BC^2D^4 + 10A^2BC^2D^3 - A^2C^{16} - 6A^2C^{13}D^2 +
        6A^2C^{13}D - 6A^2C^{13} - 9A^2C^{10}D^4 + 30A^2C^{10}D^3 -
        24A^2C^{10}D^2 + 30A^2C^{10}D - 9A^2C^{10} - 4A^2C^7D^6 +
        30A^2C^7D^5 - 30A^2C^7D^4 + 32A^2C^7D^3 - 30A^2C^7D^2 +
        30A^2C^7D - 4A^2C^7 + 15A^2C^4D^6 - 70A^2C^4D^5 +
        85A^2C^4D^4 - 70A^2C^4D^3 + 15A^2C^4D^2 + 6A^2CD^5 +
        AB^6C^5 - AB^5C^7D^3 - 3AB^5C^7D^2 + 9AB^5C^7D +
        5AB^5C^7 - 5AB^5C^4D + 12AB^4C^9D^2 + 34AB^4C^9D +
        AB^4C^9 + AB^4C^6D^4 + 8AB^4C^6D^3 - 53AB^4C^6D^2 -
        2AB^4C^6D - 9AB^4C^6 + 10AB^4C^3D^2 + AB^3C^{11}D^3 +
        15AB^3C^{11}D^2 + 9AB^3C^{11}D - 21AB^3C^{11} + 3AB^3C^8D^4 -
        58AB^3C^8D^3 - 101AB^3C^8D^2 + 23AB^3C^8D - 37AB^3C^8 -
        6AB^3C^5D^4 + 123AB^3C^5D^3 - 39AB^3C^5D^2 + 48AB^3C^5D
        - 6AB^3C^5 - 10AB^3C^2D^3 - 18AB^2C^{13}D - 15AB^2C^{13} -
        5AB^2C^{10}D^4 - 49AB^2C^{10}D^3 + 39AB^2C^{10}D^2 +
        36AB^2C^{10}D - 11AB^2C^{10} - AB^2C^7D^6 - 6AB^2C^7D^5 +
        115AB^2C^7D^4 + 52AB^2C^7D^3 - 15AB^2C^7D^2 + 74AB^2C^7D
        + 9AB^2C^7 - 135AB^2C^4D^4 + 80AB^2C^4D^3 - 90AB^2C^4D^2
        + 15AB^2C^4D + 5AB^2CD^4 - 2ABC^{15}D + 4ABC^{15} +
        ABC^{12}D^3 + 44ABC^{12}D^2 - 9ABC^{12}D + 27ABC^{12} +
        7ABC^9D^5 + 51ABC^9D^4 - 109ABC^9D^3 + 57ABC^9D^2 -
        69ABC^9D + 31ABC^9 + ABC^6D^7 + 8ABC^6D^6 -
        109ABC^6D^5 + 69ABC^6D^4 - 79ABC^6D^3 + 25ABC^6D^2 -
        56ABC^6D + 7ABC^6 + 75ABC^3D^5 - 65ABC^3D^4 +
        70ABC^3D^3 - 10ABC^3D^2 - ABD^5 + AC^{17} + 3AC^{14}D^2 -
        9AC^{14}D + 3AC^{14} - 30AC^{11}D^3 + 30AC^{11}D^2 - 30AC^{11}D -
        2AC^8D^6 - 24AC^8D^5 + 78AC^8D^4 - 86AC^8D^3 + 78AC^8D^2
        - 24AC^8D - 2AC^8 - 6AC^5D^7 + 48AC^5D^6 - 76AC^5D^5 +
        97AC^5D^4 - 76AC^5D^3 + 48AC^5D^2 - 6AC^5D - 15AC^2D^6 +
        15AC^2D^5 - 15AC^2D^4 - B^6C^6D^4 - B^6C^6D^3 - B^6C^6D^2 -
        B^6C^6D - B^6C^6 - 2B^5C^8D^3 - 6B^5C^8D^2 - 10B^5C^8D -
        2B^5C^8 + 6B^5C^5D^4 + 6B^5C^5D^3 + 6B^5C^5D^2 + 5B^5C^5D
        + B^5C^5 - 3B^4C^{10}D^2 - 3B^4C^{10}D + 9B^4C^{10} + 3B^4C^7D^5 +
        6B^4C^7D^4 + 19B^4C^7D^3 + 38B^4C^7D^2 + 7B^4C^7D +
        12B^4C^7 - 15B^4C^4D^4 - 15B^4C^4D^3 - 10B^4C^4D^2 -
        5B^4C^4D + B^3C^{12}D^2 + 12B^3C^{12}D + 12B^3C^{12} + 5B^3C^9D^4
        + 14B^3C^9D^3 - 27B^3C^9D^2 - 23B^3C^9D + 4B^3C^9 -
        B^3C^6D^7 - B^3C^6D^6 - 18B^3C^6D^5 - 41B^3C^6D^4 -
        14B^3C^6D^3 - 27B^3C^6D^2 - 31B^3C^6D - 7B^3C^6 +
        20B^3C^3D^4 + 10B^3C^3D^3 + 10B^3C^3D^2 + B^2C^{14}D -
        3B^2C^{14} - B^2C^{11}D^4 - 13B^2C^{11}D^3 - 27B^2C^{11}D^2 -
        3B^2C^{11}D - 28B^2C^{11} - 2B^2C^8D^6 - 15B^2C^8D^5 +
        13B^2C^8D^4 + 39B^2C^8D^3 - 6B^2C^8D^2 + 59B^2C^8D -
        25B^2C^8 + 6B^2C^5D^7 + 6B^2C^5D^6 + 51B^2C^5D^5 -
        24B^2C^5D^4 + 43B^2C^5D^3 - 9B^2C^5D^2 + 40B^2C^5D -
        3B^2C^5 - 5B^2C^2D^4 - 10B^2C^2D^3 - BC^{16} - BC^{13}D^3 +
        BC^{13}D^2 + 5BC^{13}D + 2BC^{10}D^5 + 22BC^{10}D^4 + 7BC^{10}D^3 +
        9BC^{10}D^2 + 15BC^{10}D + 14BC^{10} + 3BC^7D^7 + 13BC^7D^6 -
        29BC^7D^5 - 28BC^7D^4 + 33BC^7D^3 - 58BC^7D^2 - 2BC^7D +
        11BC^7 - 15BC^4D^7 - 15BC^4D^6 + 60BC^4D^5 - 120BC^4D^4 +
        110BC^4D^3 - 65BC^4D^2 + 5BC^4D - 5BCD^5 + 5BCD^4 +
        C^{15}D^2 - C^{15}D + C^{15} + 3C^{12}D^4 - 12C^{12}D^3 + 15C^{12}D^2 -
        12C^{12}D + 3C^{12} + 2C^9D^6 - 30C^9D^5 + 60C^9D^4 - 86C^9D^3 +
        60C^9D^2 - 30C^9D + 2C^9 + C^6D^8 - 22C^6D^7 + 79C^6D^6 -
        134C^6D^5 + 170C^6D^4 - 134C^6D^3 + 79C^6D^2 - 22C^6D + C^6 +
        20C^3D^7 - 55C^3D^6 + 75C^3D^5 - 55C^3D^4 + 20C^3D^3 + D^6)u^3
+( A^2B^6C^6D^3 + 2A^2B^6C^6D^2 + 3A^2B^6C^6D - A^2B^6C^6 +
        2A^2B^5C^8D^2 - A^2B^5C^8D - 7A^2B^5C^8 + A^2B^5C^5D^4 -
        5A^2B^5C^5D^3 - 11A^2B^5C^5D^2 + 7A^2B^5C^5D - 3A^2B^5C^5
        - A^2B^4C^{10}D^2 - 10A^2B^4C^{10}D - 7A^2B^4C^{10} -
        3A^2B^4C^7D^4 - 6A^2B^4C^7D^3 + 19A^2B^4C^7D^2 +
        25A^2B^4C^7D - 2A^2B^4C^7 - 5A^2B^4C^4D^4 +
        10A^2B^4C^4D^3 - 20A^2B^4C^4D^2 + 10A^2B^4C^4D -
        A^2B^3C^{12}D + 3A^2B^3C^{12} - 5A^2B^3C^9D^3 + 9A^2B^3C^9D^2
        + 2A^2B^3C^9 + A^2B^3C^6D^6 - A^2B^3C^6D^5 + 13A^2B^3C^6D^4
        - 12A^2B^3C^6D^3 - 44A^2B^3C^6D^2 - 4A^2B^3C^6D -
        A^2B^3C^6 + 10A^2B^3C^3D^4 + 30A^2B^3C^3D^3 -
        10A^2B^3C^3D^2 + A^2B^2C^{14} - A^2B^2C^{11}D^3 +
        10A^2B^2C^{11}D^2 - 10A^2B^2C^{11}D - 7A^2B^2C^{11} +
        2A^2B^2C^8D^5 + 25A^2B^2C^8D^4 - 16A^2B^2C^8D^3 +
        43A^2B^2C^8D^2 + 29A^2B^2C^8D - 5A^2B^2C^8 + A^2B^2C^5D^7
        - 5A^2B^2C^5D^6 + 4A^2B^2C^5D^5 - 23A^2B^2C^5D^4 +
        21A^2B^2C^5D^3 - 14A^2B^2C^5D^2 + 24A^2B^2C^5D -
        2A^2B^2C^5 - 25A^2B^2C^2D^4 + 4A^2BC^{13}D^2 + A^2BC^{13}D +
        8A^2BC^{13} + 6A^2BC^{10}D^4 - 25A^2BC^{10}D^3 + 16A^2BC^{10}D^2
        - 25A^2BC^{10}D + 26A^2BC^{10} - 37A^2BC^7D^5 + 19A^2BC^7D^4
        - 31A^2BC^7D^3 - 21A^2BC^7D^2 - 47A^2BC^7D + 15A^2BC^7 -
        5A^2BC^4D^7 + 10A^2BC^4D^6 + 40A^2BC^4D^5 - 80A^2BC^4D^4
        + 120A^2BC^4D^3 - 60A^2BC^4D^2 + 5A^2BC^4D + 5A^2BCD^4 -
        2A^2C^{15}D + A^2C^{15} - 9A^2C^{12}D^3 + 15A^2C^{12}D^2 -
        20A^2C^{12}D + A^2C^{12} - 11A^2C^9D^5 + 51A^2C^9D^4 -
        81A^2C^9D^3 + 89A^2C^9D^2 - 39A^2C^9D - A^2C^9 - 5A^2C^6D^7
        + 46A^2C^6D^6 - 103A^2C^6D^5 + 171A^2C^6D^4 - 167A^2C^6D^3 +
        121A^2C^6D^2 - 28A^2C^6D + A^2C^6 + 10A^2C^3D^7 -
        50A^2C^3D^6 + 75A^2C^3D^5 - 70A^2C^3D^4 + 20A^2C^3D^3 +
        A^2D^6 + AB^7C^5 - AB^6C^7D^3 - 3AB^6C^7D^2 + 4AB^6C^7D +
        5AB^6C^7 - 5AB^6C^4D + 3AB^5C^9D^2 + 18AB^5C^9D +
        6AB^5C^9 + 6AB^5C^6D^3 - 26AB^5C^6D^2 - 12AB^5C^6D -
        3AB^5C^6 + 10AB^5C^3D^2 + 3AB^4C^{11}D^2 + 9AB^4C^{11}D -
        5AB^4C^{11} + 3AB^4C^8D^4 - 3AB^4C^8D^3 - 57AB^4C^8D^2 -
        12AB^4C^8 + AB^4C^5D^5 + 60AB^4C^5D^3 + 9AB^4C^5D^2 +
        17AB^4C^5D - 2AB^4C^5 - 10AB^4C^2D^3 - 2AB^3C^{13}D -
        6AB^3C^{13} + AB^3C^{10}D^4 - 12AB^3C^{10}D^2 + 36AB^3C^{10}D -
        3AB^3C^{10} - AB^3C^7D^6 - 3AB^3C^7D^5 - 19AB^3C^7D^4 +
        57AB^3C^7D^3 - 48AB^3C^7D^2 + 31AB^3C^7D - 2AB^3C^7 -
        5AB^3C^4D^5 - 50AB^3C^4D^4 - 10AB^3C^4D^3 - 25AB^3C^4D^2
        + 5AB^3CD^4 - 3AB^2C^{12}D^3 - 16AB^2C^{12}D^2 + 6AB^2C^{12}D
        - 9AB^2C^{12} - 7AB^2C^9D^5 - 21AB^2C^9D^4 + 57AB^2C^9D^3 -
        91AB^2C^9D^2 + 18AB^2C^9D - 27AB^2C^9 + 3AB^2C^6D^6 +
        64AB^2C^6D^5 - 90AB^2C^6D^4 + 147AB^2C^6D^3 -
        49AB^2C^6D^2 + 33AB^2C^6D - 9AB^2C^6 + 10AB^2C^3D^5 +
        15AB^2C^3D^4 + 10AB^2C^3D^2 - AB^2D^5 + 4ABC^{14}D -
        9ABC^{14} + 11ABC^{11}D^4 + 28ABC^{11}D^3 - 23ABC^{11}D^2 +
        69ABC^{11}D - 15ABC^{11} + 13ABC^8D^6 + 11ABC^8D^5 -
        77ABC^8D^4 + 139ABC^8D^3 - 137ABC^8D^2 + 101ABC^8D -
        2ABC^8 + ABC^5D^8 - 74ABC^5D^6 + 129ABC^5D^5 -
        211ABC^5D^4 + 147ABC^5D^3 - 124ABC^5D^2 + 31ABC^5D -
        ABC^5 + 15ABC^2D^6 - 25ABC^2D^5 + 25ABC^2D^4 -
        10ABC^2D^3 + AC^{16}D + AC^{16} - 2AC^{13}D^2 - 3AC^{13}D +
        8AC^{13} - 7AC^{10}D^5 - 9AC^{10}D^4 - 5AC^{10}D^2 - 31AC^{10}D +
        12AC^{10} - 5AC^7D^7 - 6AC^7D^6 + 48AC^7D^5 - 76AC^7D^4 +
        62AC^7D^3 - 6AC^7D^2 - 27AC^7D + 4AC^7 - 5AC^4D^8 +
        45AC^4D^7 - 100AC^4D^6 + 155AC^4D^5 - 135AC^4D^4 +
        85AC^4D^3 - 15AC^4D^2 - 5ACD^7 + 6ACD^6 - 6ACD^5 -
        B^7C^6D^4 - B^7C^6D^3 - B^7C^6D^2 - B^7C^6D - B^7C^6 -
        B^6C^8D^3 - 4B^6C^8D^2 - 7B^6C^8D - 3B^6C^8 - B^6C^5D^5 +
        6B^6C^5D^4 + 6B^6C^5D^3 + 6B^6C^5D^2 + 5B^6C^5D + B^6C^5 -
        B^5C^{10}D^2 - 4B^5C^{10}D + 2B^5C^{10} + 3B^5C^7D^5 + 3B^5C^7D^4
        + 7B^5C^7D^3 + 25B^5C^7D^2 + 11B^5C^7D + 7B^5C^7 +
        5B^5C^4D^5 - 15B^5C^4D^4 - 15B^5C^4D^3 - 10B^5C^4D^2 -
        5B^5C^4D + 2B^4C^{12}D + 5B^4C^{12} + 2B^4C^9D^4 + 9B^4C^9D^3
        - 5B^4C^9D^2 - 16B^4C^9D + B^4C^9 - B^4C^6D^7 + 2B^4C^6D^6 -
        14B^4C^6D^5 - 25B^4C^6D^4 + 10B^4C^6D^3 - 21B^4C^6D^2 -
        17B^4C^6D - 4B^4C^6 - 10B^4C^3D^5 + 20B^4C^3D^4 +
        10B^4C^3D^3 + 10B^4C^3D^2 + 8B^3C^{11}D^2 + 2B^3C^{11}D -
        2B^3C^{11} - B^3C^8D^6 - 5B^3C^8D^5 + 11B^3C^8D^4 -
        20B^3C^8D^3 + 17B^3C^8D^2 + 34B^3C^8D + 3B^3C^8 - B^3C^5D^8
        + 6B^3C^5D^7 - 8B^3C^5D^6 + 27B^3C^5D^5 + 7B^3C^5D^4 -
        20B^3C^5D^3 - 10B^3C^5D^2 + 19B^3C^5D + 10B^3C^2D^5 -
        5B^3C^2D^4 - 10B^3C^2D^3 - B^2C^{13}D^2 - 3B^2C^{13}D +
        8B^2C^{13} - 2B^2C^{10}D^5 - 17B^2C^{10}D^4 - 5B^2C^{10}D^3 -
        6B^2C^{10}D^2 - 48B^2C^{10}D + 11B^2C^{10} - 11B^2C^7D^6 +
        45B^2C^7D^5 - 64B^2C^7D^4 + 46B^2C^7D^3 + 8B^2C^7D^2 -
        57B^2C^7D + 5B^2C^4D^8 - 15B^2C^4D^7 + 15B^2C^4D^6 +
        20B^2C^4D^5 - 65B^2C^4D^4 + 80B^2C^4D^3 - 10B^2C^4D^2 -
        5B^2C^4D - 10B^2CD^5 + 5B^2CD^4 - BC^{15}D - BC^{15} +
        BC^{12}D^4 + BC^{12}D^3 - 5BC^{12}D^2 + 11BC^{12}D - 14BC^{12} +
        8BC^9D^6 + 8BC^9D^5 + 6BC^9D^4 + 46BC^9D^3 - 40BC^9D^2 +
        74BC^9D - 20BC^9 + 5BC^6D^8 - 9BC^6D^7 + 11BC^6D^6 -
        41BC^6D^5 + 37BC^6D^4 - 7BC^6D^3 - 44BC^6D^2 + 54BC^6D -
        6BC^6 - 10BC^3D^8 + 20BC^3D^7 - 20BC^3D^6 - 15BC^3D^5 +
        45BC^3D^4 - 50BC^3D^3 + 10BC^3D^2 + BD^6 + BD^5 + C^{14}D^3 +
        C^{14} + C^{11}D^5 - 6C^{11}D^4 + 6C^{11}D^3 + C^{11}D^2 - 6C^{11}D + 6C^{11}
        - C^8D^7 - 21C^8D^6 + 51C^8D^5 - 87C^8D^4 + 57C^8D^3 -
        21C^8D^2 - 14C^8D + 6C^8 + C^5D^9 - 24C^5D^8 + 105C^5D^7 -
        231C^5D^6 + 339C^5D^5 - 321C^5D^4 + 200C^5D^3 - 69C^5D^2 +
        6C^5D + 10C^2D^8 - 30C^2D^7 + 45C^2D^6 - 36C^2D^5 +
        15C^2D^4)u^2
+ (A^2B^6C^5D^4 + A^2B^6C^5D^3 + A^2B^6C^5D^2 + A^2B^6C^5D +
        2A^2B^6C^5 + 2A^2B^5C^7D^3 + 5A^2B^5C^7D^2 + 17A^2B^5C^7D
        + 6A^2B^5C^7 - 5A^2B^5C^4D^4 - 5A^2B^5C^4D^3 -
        5A^2B^5C^4D^2 - 10A^2B^5C^4D + 8A^2B^4C^9D^2 +
        20A^2B^4C^9D - 10A^2B^4C^9 - 3A^2B^4C^6D^5 -
        6A^2B^4C^6D^4 - 15A^2B^4C^6D^3 - 80A^2B^4C^6D^2 -
        2A^2B^4C^6D - 19A^2B^4C^6 + 10A^2B^4C^3D^4 +
        10A^2B^4C^3D^3 + 20A^2B^4C^3D^2 + A^2B^3C^{11}D^2 -
        13A^2B^3C^{11}D - 24A^2B^3C^{11} - 5A^2B^3C^8D^4 -
        41A^2B^3C^8D^3 - 30A^2B^3C^8D^2 + 47A^2B^3C^8D -
        40A^2B^3C^8 + A^2B^3C^5D^7 + A^2B^3C^5D^6 + 15A^2B^3C^5D^5
        + 33A^2B^3C^5D^4 + 132A^2B^3C^5D^3 - 45A^2B^3C^5D^2 +
        87A^2B^3C^5D - 9A^2B^3C^5 - 10A^2B^3C^2D^4 -
        20A^2B^3C^2D^3 - 4A^2B^2C^{13}D + A^2B^2C^{10}D^4 +
        4A^2B^2C^{10}D^3 + 49A^2B^2C^{10}D^2 + 23A^2B^2C^{10}D +
        18A^2B^2C^{10} + 2A^2B^2C^7D^6 + 14A^2B^2C^7D^5 +
        48A^2B^2C^7D^4 - 8A^2B^2C^7D^3 - 17A^2B^2C^7D^2 +
        39A^2B^2C^7D + 24A^2B^2C^7 - 5A^2B^2C^4D^7 -
        5A^2B^2C^4D^6 - 35A^2B^2C^4D^5 - 130A^2B^2C^4D^4 +
        100A^2B^2C^4D^3 - 140A^2B^2C^4D^2 + 20A^2B^2C^4D +
        10A^2B^2CD^4 + 2A^2BC^{15} + A^2BC^{12}D^3 + 9A^2BC^{12}D^2 -
        10A^2BC^{12}D + 14A^2BC^{12} - 2A^2BC^9D^5 - 7A^2BC^9D^4 -
        62A^2BC^9D^3 + 31A^2BC^9D^2 - 70A^2BC^9D + 21A^2BC^9 -
        3A^2BC^6D^7 - 9A^2BC^6D^6 - 38A^2BC^6D^5 + 77A^2BC^6D^4
        - 96A^2BC^6D^3 + 67A^2BC^6D^2 - 75A^2BC^6D + 8A^2BC^6 +
        10A^2BC^3D^7 + 10A^2BC^3D^6 + 50A^2BC^3D^5 -
        45A^2BC^3D^4 + 80A^2BC^3D^3 - 10A^2BC^3D^2 - A^2BD^5 -
        A^2C^{14}D^2 - A^2C^{14}D - 3A^2C^{11}D^4 + 3A^2C^{11}D^3 -
        8A^2C^{11}D - 2A^2C^{11} - 2A^2C^8D^6 + 19A^2C^8D^5 -
        9A^2C^8D^4 + 5A^2C^8D^3 + 29A^2C^8D^2 - 9A^2C^8D - 3A^2C^8
        - A^2C^5D^8 + 17A^2C^5D^7 - 33A^2C^5D^6 + 31A^2C^5D^5 +
        A^2C^5D^4 - 33A^2C^5D^3 + 42A^2C^5D^2 - 6A^2C^5D -
        10A^2C^2D^7 + 5A^2C^2D^6 - 15A^2C^2D^4 - AB^6C^6D^4 -
        2AB^6C^6D^3 - 3AB^6C^6D^2 - 9AB^6C^6D + AB^6C^6 -
        2AB^5C^8D^3 - 18AB^5C^8D^2 - 18AB^5C^8D + 14AB^5C^8 +
        AB^5C^5D^5 + 6AB^5C^5D^4 + 12AB^5C^5D^3 + 42AB^5C^5D^2 -
        12AB^5C^5D + 9AB^5C^5 - AB^4C^{10}D^3 - 9AB^4C^{10}D^2 +
        24AB^4C^{10}D + 31AB^4C^{10} + 3AB^4C^7D^5 + 12AB^4C^7D^4 +
        73AB^4C^7D^3 + 21AB^4C^7D^2 - 33AB^4C^7D + 21AB^4C^7 -
        5AB^4C^4D^5 - 15AB^4C^4D^4 - 75AB^4C^4D^3 + 40AB^4C^4D^2
        - 35AB^4C^4D + 3AB^3C^{12}D^2 + 28AB^3C^{12}D + 6AB^3C^{12} +
        5AB^3C^9D^4 + 20AB^3C^9D^3 - 129AB^3C^9D^2 - 20AB^3C^9D
        - 32AB^3C^9 - AB^3C^6D^7 - 5AB^3C^6D^6 - 24AB^3C^6D^5 -
        128AB^3C^6D^4 + 70AB^3C^6D^3 - 42AB^3C^6D^2 + AB^3C^6D -
        28AB^3C^6 + 10AB^3C^3D^5 + 60AB^3C^3D^4 - 60AB^3C^3D^3 +
        50AB^3C^3D^2 - 11AB^2C^{14} - 3AB^2C^{11}D^4 - 23AB^2C^{11}D^3
        - 62AB^2C^{11}D^2 + 27AB^2C^{11}D - 52AB^2C^{11} - 2AB^2C^8D^6
        - 10AB^2C^8D^5 + 21AB^2C^8D^4 + 184AB^2C^8D^3 -
        92AB^2C^8D^2 + 153AB^2C^8D - 41AB^2C^8 + AB^2C^5D^8 +
        6AB^2C^5D^7 + 27AB^2C^5D^6 + 87AB^2C^5D^5 -
        135AB^2C^5D^4 + 173AB^2C^5D^3 - 115AB^2C^5D^2 +
        90AB^2C^5D - 5AB^2C^5 - 25AB^2C^2D^5 + 45AB^2C^2D^4 -
        30AB^2C^2D^3 - ABC^{16} + 2ABC^{13}D^3 + 5ABC^{13}D^2 +
        11ABC^{13}D - 2ABC^{13} + 10ABC^{10}D^5 + 22ABC^{10}D^4 +
        25ABC^{10}D^3 + ABC^{10}D^2 + 51ABC^{10}D + 11ABC^{10} +
        6ABC^7D^7 - 7ABC^7D^6 - 32ABC^7D^5 - 76ABC^7D^4 +
        44ABC^7D^3 - 112ABC^7D^2 + 28ABC^7D + 11ABC^7 -
        5ABC^4D^8 - 15ABC^4D^7 - 15ABC^4D^6 + 70ABC^4D^5 -
        125ABC^4D^4 + 110ABC^4D^3 - 80ABC^4D^2 + 5ABC^4D -
        11ABCD^5 + 5ABCD^4 + AC^{15}D + AC^{15} - 3AC^{12}D^4 -
        5AC^{12}D^2 - 3AC^{12}D + 5AC^{12} - 7AC^9D^6 + 2AC^9D^5 -
        12AC^9D^4 - 2AC^9D^3 - 8AC^9D^2 - 22AC^9D + 5AC^9 -
        3AC^6D^8 + 7AC^6D^7 + 9AC^6D^6 + 6AC^6D^5 + 2AC^6D^4 +
        20AC^6D^3 + 9AC^6D^2 - 15AC^6D + AC^6 + 10AC^3D^8 -
        20AC^3D^7 + 15AC^3D^6 + 5AC^3D^5 - 20AC^3D^4 + 20AC^3D^3
        + AD^7 + AD^6 - B^7C^5D^5 - B^6C^7D^4 - B^6C^7D^3 +
        2B^6C^7D^2 + 3B^6C^7D - 3B^6C^7 + 5B^6C^4D^5 - B^5C^9D^3 -
        4B^5C^9D^2 - 13B^5C^9D - 15B^5C^9 + 3B^5C^6D^6 +
        2B^5C^6D^5 + B^5C^6D^4 + 4B^5C^6D^3 - 4B^5C^6D^2 +
        10B^5C^6D - 4B^5C^6 - 10B^5C^3D^5 - 2B^4C^{11}D^2 -
        16B^4C^{11}D - 7B^4C^{11} + 2B^4C^8D^5 + 9B^4C^8D^4 +
        15B^4C^8D^3 + 68B^4C^8D^2 + 37B^4C^8D + 19B^4C^8 -
        B^4C^5D^8 - 13B^4C^5D^6 - 12B^4C^5D^5 + 4B^4C^5D^4 -
        38B^4C^5D^3 - 9B^4C^5D^2 + 2B^4C^5D + 7B^4C^5 +
        10B^4C^2D^5 + B^3C^{13}D + 9B^3C^{13} + 3B^3C^{10}D^4 +
        29B^3C^{10}D^3 + 32B^3C^{10}D^2 - 3B^3C^{10}D + 40B^3C^{10} -
        B^3C^7D^7 - 7B^3C^7D^6 + B^3C^7D^5 - 90B^3C^7D^4 -
        76B^3C^7D^3 - 31B^3C^7D^2 - 92B^3C^7D + 14B^3C^7 +
        5B^3C^4D^8 + 25B^3C^4D^6 + 30B^3C^4D^5 + 45B^3C^4D^4 +
        5B^3C^4D^3 + 35B^3C^4D^2 - 25B^3C^4D - 5B^3CD^5 + B^2C^{15} -
        3B^2C^{12}D^3 - 12B^2C^{12}D^2 - 4B^2C^{12}D - 6B^2C^{12} -
        B^2C^9D^6 - 20B^2C^9D^5 - 35B^2C^9D^4 + 8B^2C^9D^3 -
        68B^2C^9D^2 - 7B^2C^9D - 34B^2C^9 + 2B^2C^6D^8 -
        5B^2C^6D^7 + 40B^2C^6D^6 + 64B^2C^6D^5 + 28B^2C^6D^4 +
        99B^2C^6D^3 + 26B^2C^6D^2 + 39B^2C^6D - 17B^2C^6 -
        10B^2C^3D^8 - 30B^2C^3D^6 - 80B^2C^3D^5 + 50B^2C^3D^4 -
        80B^2C^3D^3 + 30B^2C^3D^2 + B^2D^5 - BC^{14}D^2 - 2BC^{14} +
        2BC^{11}D^5 + 4BC^{11}D^4 + 11BC^{11}D^3 - 6BC^{11}D^2 + 17BC^{11}D
        - 7BC^{11} + 6BC^8D^7 + 7BC^8D^6 + 28BC^8D^5 - 33BC^8D^4 +
        96BC^8D^3 - 53BC^8D^2 + 62BC^8D - 2BC^8 + 2BC^5D^9 -
        12BC^5D^8 + 19BC^5D^7 - 128BC^5D^6 + 156BC^5D^5 -
        241BC^5D^4 + 146BC^5D^3 - 119BC^5D^2 + 37BC^5D - BC^5 +
        10BC^2D^8 + 40BC^2D^6 - 20BC^2D^5 + 40BC^2D^4 - 10BC^2D^3
        + C^{13}D^3 + C^{13} - 2C^{10}D^6 + C^{10}D^5 - 6C^{10}D^4 + 2C^{10}D^3 +
        C^{10}D^2 - 6C^{10}D + 4C^{10} - 3C^7D^8 - C^7D^7 - 3C^7D^5 -
        15C^7D^4 + 3C^7D^3 - 14C^7D + 3C^7 - 5C^4D^9 + 30C^4D^8 -
        45C^4D^7 + 45C^4D^6 + 15C^4D^5 - 45C^4D^4 + 50C^4D^3 -
        15C^4D^2 - 5CD^8 - CD^6 - 6CD^5)u
 -A^2B^5C^6D^3 - 2A^2B^5C^6D^2 - 3A^2B^5C^6D + 2A^2B^5C^6 - 3A^2B^4C^8D^2 + 8A^2B^4C^8D + 11A^2B^4C^8 + 6A^2B^4C^5D^3 + 12A^2B^4C^5D^2 - 12A^2B^4C^5D + 4A^2B^4C^5 +
        6A^2B^3C^{10}D^2 + 27A^2B^3C^{10}D + 6A^2B^3C^{10} + 3A^2B^3C^7D^4 + 10A^2B^3C^7D^3 - 61A^2B^3C^7D^2 - 20A^2B^3C^7D - 5A^2B^3C^7 - 15A^2B^3C^4D^3 + 30A^2B^3C^4D^2 -
        15A^2B^3C^4D + 2A^2B^2C^{12}D^2 - 15A^2B^2C^{12} - 22A^2B^2C^9D^3 - 66A^2B^2C^9D^2 + 24A^2B^2C^9D - 38A^2B^2C^9 - A^2B^2C^6D^6 - 2A^2B^2C^6D^5 - 21A^2B^2C^6D^4 +
        130A^2B^2C^6D^3 - 28A^2B^2C^6D^2 + 60A^2B^2C^6D - 15A^2B^2C^6 - 40A^2B^2C^3D^3 + 20A^2B^2C^3D^2 - 3A^2BC^{14}D - 4A^2BC^{14} - 8A^2BC^{11}D^3 + 12A^2BC^{11}D^2 + 30A^2BC^{11}D -
        3A^2BC^{11} - 3A^2BC^8D^5 + 36A^2BC^8D^4 + 47A^2BC^8D^3 - 66A^2BC^8D^2 + 69A^2BC^8D + 4A^2BC^8 + 6A^2BC^5D^6 + 12A^2BC^5D^5 - 131A^2BC^5D^4 + 122A^2BC^5D^3 -
        135A^2BC^5D^2 + 36A^2BC^5D - A^2BC^5 + 30A^2BC^2D^4 - 10A^2BC^2D^3 + A^2C^{16} + 6A^2C^{13}D^2 - 6A^2C^{13}D + 6A^2C^{13} + 9A^2C^{10}D^4 - 30A^2C^{10}D^3 + 24A^2C^{10}D^2 - 30A^2C^{10}D +
        9A^2C^{10} + 4A^2C^7D^6 - 30A^2C^7D^5 + 30A^2C^7D^4 - 32A^2C^7D^3 + 30A^2C^7D^2 - 30A^2C^7D + 4A^2C^7 - 15A^2C^4D^6 + 70A^2C^4D^5 - 85A^2C^4D^4 + 70A^2C^4D^3 - 15A^2C^4D^2 -
        6A^2CD^5 - AB^6C^5 + AB^5C^7D^3 + 3AB^5C^7D^2 - 9AB^5C^7D - 5AB^5C^7 + 5AB^5C^4D - 12AB^4C^9D^2 - 34AB^4C^9D - AB^4C^9 - AB^4C^6D^4 - 8AB^4C^6D^3 + 53AB^4C^6D^2 +
        2AB^4C^6D + 9AB^4C^6 - 10AB^4C^3D^2 - AB^3C^{^11}D^3 - 15AB^3C^{11}D^2 - 9AB^3C^{11}D + 21AB^3C^{11} - 3AB^3C^8D^4 + 58AB^3C^8D^3 + 101AB^3C^8D^2 - 23AB^3C^8D + 37AB^3C^8 +
        6AB^3C^5D^4 - 123AB^3C^5D^3 + 39AB^3C^5D^2 - 48AB^3C^5D + 6AB^3C^5 + 10AB^3C^2D^3 + 18AB^2C^{13}D + 15AB^2C^{13} + 5AB^2C^{10}D^4 + 49AB^2C^{10}D^3 - 39AB^2C^{10}D^2 -
        36AB^2C^{10}D + 11AB^2C^{10} + AB^2C^7D^6 + 6AB^2C^7D^5 - 115AB^2C^7D^4 - 52AB^2C^7D^3 + 15AB^2C^7D^2 - 74AB^2C^7D - 9AB^2C^7 + 135AB^2C^4D^4 - 80AB^2C^4D^3 + 90AB^2C^4D^2
        - 15AB^2C^4D - 5AB^2CD^4 + 2ABC^{15}D - 4ABC^{15} - ABC^{12}D^3 - 44ABC^{12}D^2 + 9ABC^{12}D - 27ABC^{12} - 7ABC^9D^5 - 51ABC^9D^4 + 109ABC^9D^3 - 57ABC^9D^2 + 69ABC^9D -
        31ABC^9 - ABC^6D^7 - 8ABC^6D^6 + 109ABC^6D^5 - 69ABC^6D^4 + 79ABC^6D^3 - 25ABC^6D^2 + 56ABC^6D - 7ABC^6 - 75ABC^3D^5 + 65ABC^3D^4 - 70ABC^3D^3 + 10ABC^3D^2 +
        ABD^5 - AC^{17} - 3AC^{14}D^2 + 9AC^{14}D - 3AC^{14} + 30AC^{11}D^3 - 30AC^{11}D^2 + 30AC^{11}D + 2AC^8D^6 + 24AC^8D^5 - 78AC^8D^4 + 86AC^8D^3 - 78AC^8D^2 + 24AC^8D + 2AC^8 + 6AC^5D^7
        - 48AC^5D^6 + 76AC^5D^5 - 97AC^5D^4 + 76AC^5D^3 - 48AC^5D^2 + 6AC^5D + 15AC^2D^6 - 15AC^2D^5 + 15AC^2D^4 + B^6C^6D^4 + B^6C^6D^3 + B^6C^6D^2 + B^6C^6D + B^6C^6 + 2B^5C^8D^3
        + 6B^5C^8D^2 + 10B^5C^8D + 2B^5C^8 - 6B^5C^5D^4 - 6B^5C^5D^3 - 6B^5C^5D^2 - 5B^5C^5D - B^5C^5 + 3B^4C^{10}D^2 + 3B^4C^{10}D - 9B^4C^{10} - 3B^4C^7D^5 - 6B^4C^7D^4 - 19B^4C^7D^3 -
        38B^4C^7D^2 - 7B^4C^7D - 12B^4C^7 + 15B^4C^4D^4 + 15B^4C^4D^3 + 10B^4C^4D^2 + 5B^4C^4D - B^3C^{12}D^2 - 12B^3C^{12}D - 12B^3C^{12} - 5B^3C^9D^4 - 14B^3C^9D^3 + 27B^3C^9D^2 +
        23B^3C^9D - 4B^3C^9 + B^3C^6D^7 + B^3C^6D^6 + 18B^3C^6D^5 + 41B^3C^6D^4 + 14B^3C^6D^3 + 27B^3C^6D^2 + 31B^3C^6D + 7B^3C^6 - 20B^3C^3D^4 - 10B^3C^3D^3 - 10B^3C^3D^2 - B^2C^{14}D
        + 3B^2C^{14} + B^2C^{11}D^4 + 13B^2C^{11}D^3 + 27B^2C^{11}D^2 + 3B^2C^{11}D + 28B^2C^{11} + 2B^2C^8D^6 + 15B^2C^8D^5 - 13B^2C^8D^4 - 39B^2C^8D^3 + 6B^2C^8D^2 - 59B^2C^8D + 25B^2C^8 - 6B^2C^5D^7 - 6B^2C^5D^6 - 51B^2C^5D^5 + 24B^2C^5D^4 - 43B^2C^5D^3 + 9B^2C^5D^2 - 40B^2C^5D + 3B^2C^5 + 5B^2C^2D^4 + 10B^2C^2D^3 + BC^{16} + BC^{13}D^3 - BC^{13}D^2 - 5BC^{13}D -
        2BC^{10}D^5 - 22BC^{10}D^4 - 7BC^{10}D^3 - 9BC^{10}D^2 - 15BC^{10}D - 14BC^{10} - 3BC^7D^7 - 13BC^7D^6 + 29BC^7D^5 + 28BC^7D^4 - 33BC^7D^3 + 58BC^7D^2 + 2BC^7D - 11BC^7 + 15BC^4D^7
        + 15BC^4D^6 - 60BC^4D^5 + 120BC^4D^4 - 110BC^4D^3 + 65BC^4D^2 - 5BC^4D + 5BCD^5 - 5BCD^4 - C^{15}D^2 + C^{15}D - C^{15} - 3C^{12}D^4 + 12C^{12}D^3 - 15C^{12}D^2 + 12C^{12}D - 3C^{12} -
        2C^9D^6 + 30C^9D^5 - 60C^9D^4 + 86C^9D^3 - 60C^9D^2 + 30C^9D - 2C^9 - C^6D^8 + 22C^6D^7 - 79C^6D^6 + 134C^6D^5 - 170C^6D^4 + 134C^6D^3 - 79C^6D^2 + 22C^6D - C^6 - 20C^3D^7 +
        55C^3D^6 - 75C^3D^5 + 55C^3D^4 - 20C^3D^3 - D^6  $.
\end{scriptsize}

\end{document}